\documentclass[a4paper,11pt,reqno]{amsart}

\usepackage{rotating,pdflscape,color,amssymb,subfigure,psfrag,amsmath,eufrak,bbm, epsfig,cancel}
\definecolor{vdarkred}{rgb}{0.6,0,0.2}
\definecolor{vdarkblue}{rgb}{0,0.2,0.6}
\usepackage[pdftex, colorlinks, linkcolor=vdarkblue,citecolor=vdarkred]{hyperref}

\date{\today}

\newcommand{\dd}{\mathrm{d}}
\newcommand{\ii}{\mathrm{i}}
\newcommand*{\deq}{\mathrel{\vcenter{\baselineskip0.65ex \lineskiplimit0pt \hbox{.}\hbox{.}}}=}
\newcommand{\Lc}{\mc{L}}

\newcommand{\ps}{\f{\partial}{\partial s}}
\newcommand{\pt}{\f{\partial}{\partial t}}

\newcommand{\ninfe}{\quad\underset{n\to\infty}{\longrightarrow}\quad}

\newcommand{\ld}{\ldots}
\newcommand{\beg}{\begin}
\newcommand{\en}{\end}

\newcommand{\trm}{\textrm}

\newcommand{\bgt}{\begin{itemize}}
\newcommand{\ent}{\end{itemize}}
\newcommand{\ite}{\item}

\newcommand{\op}{\operatorname}
\newcommand{\eqre}{\eqref}
\newcommand{\re}{\ref}
\newcommand{\la}{\label}

\newcommand{\sd}{\,\cdot\,}

\newcommand{\rfl}{\rfloor}
\newcommand{\lfl}{\lfloor}
\newcommand{\si}{\sigma}

\newcommand{\Var}{\operatorname{Var}}

\newcommand{\Hc}{\mc{H}}

\newcommand{\diag}{\operatorname{diag}}

\newcommand{\ds}{\displaystyle}

\newcommand{\Tr}{\operatorname{Tr}}

\newcommand{\ninf}{\underset{n\to\infty}{\longrightarrow}}

\newcommand{\wt}{\widehat}

\newcommand{\me}{\mathrm{e}}
\newcommand{\E}{\mathbb{E}}

\newcommand{\R}{\mathbb{R}}
\newcommand{\C}{\mathbb{C}}

\newcommand{\tS}{\widetilde{\mc{S}}}
\newcommand{\ud}{\mathrm{d}}

\newcommand{\pro}{probability }
\newcommand{\Sc}{\mc{S}}

\newcommand{\f}{\frac}
\newcommand{\ff}{\frac{1}}
\newcommand{\lf}{\left}
\newcommand{\ri}{\right}

\newcommand{\st}{such that }

\newcommand{\lam}{\lambda}
\newcommand{\La}{\Lambda}
\newcommand{\ti}{\times}

\newcommand{\vfi}{\varphi}
\newcommand{\ste}{\,;\, }
\newcommand{\mc}{\mathcal }

\newcommand{\eps}{\varepsilon}

\newcommand{\hfi}{\widehat{\phi}}

\newcommand{\bck}{\backslash}
\newcommand{\al}{\alpha}

\newcommand{\Cc}{\mc{C}}
\newcommand{\ninfP}{\quad\stackrel{P}{\underset{n\to\infty}{\longrightarrow}}\quad}
\newcommand{\ninfd}{\quad\stackrel{\mathrm{dist.}}{\underset{n\to\infty}{\longrightarrow}}\quad}
\newcommand{\cvP}{\stackrel{P}{\longrightarrow}}
\newcommand{\cvd}{\stackrel{\mathrm{dist.}}{\longrightarrow}}

\newcommand{\ovl}{\overline}

\newcommand{\bbm}{\begin{bmatrix}}
\newcommand{\ebm}{\end{bmatrix}}
\newcommand{\bes}{\begin{equation*}}
\newcommand{\ees}{\end{equation*}}
\newcommand{\be}{\begin{equation}}
\newcommand{\ee}{\end{equation}}
\newcommand{\beqy}{\begin{eqnarray}}
\newcommand{\eeqy}{\end{eqnarray}}
\newcommand{\beq}{\begin{eqnarray*}}
\newcommand{\eeq}{\end{eqnarray*}}
\newcommand{\one}{\mathbbm{1}}
\newcommand{\lto}{\longrightarrow}

\newcommand{\bpr}{\beg{pr}}
\newcommand{\epr}{\end{pr}}

\newcommand{\ie}{i.e. }

\newcommand{\bpm}{\begin{pmatrix}}
\newcommand{\epm}{\end{pmatrix}}

\newcommand{\Bo}{b_n^{\circ}}

\newcommand{\DG}{\Delta\mathrm{G}}
\newcommand{\sid}{\si_{\op{d}}}

\renewcommand{\Im}{\mathfrak{Im}}
\renewcommand{\Re}{\mathfrak{Re}}

\newtheorem{Th}{Theorem}[]

\newtheorem{propo}[Th]{Proposition}

\newtheorem{lem}[Th]{Lemma}

\theoremstyle{definition}
\newtheorem{rmq}[]{Remark}

\long\def\symbolfootnote[#1]#2{\begingroup
\def\thefootnote{\fnsymbol{footnote}}\footnote[#1]{#2}\endgroup}

\makeatletter
\newcommand*{\inlineequation}[2][]{%
  \begingroup
    \refstepcounter{equation}%
    \ifx\\#1\\%
    \else
      \label{#1}%
    \fi
    \relpenalty=10000 %
    \binoppenalty=10000 %
    \ensuremath{%
      #2%
    }%
    ~ \hfill \@eqnnum
  \endgroup
}
\makeatother

\author[Florent Benaych-Georges]{Florent Benaych-Georges}\address{Florent Benaych-Georges, Universit\'e Paris Descartes,
45, rue des Saints-P\`eres
75270 Paris Cedex 06, France.}\email{florent.benaych-georges@parisdescartes.fr}

\author[Nathana\"el Enriquez]{Nathana\"el Enriquez}\address{Nathana\"el Enriquez, MODAL'X, 200 Avenue de la R\'epublique, 92001, Nanterre and LPMA,   Case courier 188, 4 Place Jussieu, 75252 Paris Cedex 05, France.} \email{nenriquez@u-paris10.fr}

\author[Alk\'eos Micha\"{\i}l]{Alk\'eos Micha\"{\i}l}\address{Alk\'eos Micha\"{\i}l, Universit\'e Paris Descartes,
45, rue des Saints-P\`eres
75270 Paris Cedex 06, France.}\email{alkeos.michail@parisdescartes.fr}

\title{Empirical spectral distribution of a matrix under perturbation}

\begin{document}
\maketitle

\begin{abstract}We provide a perturbative expansion for the empirical spectral distribution of a Hermitian  matrix with large size perturbed by a random matrix with small operator norm whose entries in the eigenvector basis of the first one are independent with a variance profile. We prove that, depending on the order of magnitude of the perturbation, several regimes can appear, called \emph{perturbative} and \emph{semi-perturbative} regimes. Depending on the regime, the leading terms of the expansion are either related to the one-dimensional Gaussian free field or to free probability theory.\end{abstract}



\section{Introduction}It is a natural and central question, in mathematics and physics, to understand how the spectral properties of an operator are altered when the operator is subject to a small perturbation. This question is at the center of \emph{Perturbation Theory} and has been studied in many different contexts. We refer the reader to Kato's book  \cite{kato} for a thorough account on this subject. 
In this text, we   provide a perturbative expansion for the empirical spectral distribution of a    Hermitian  matrix with large size  perturbed by a   random matrix with small operator norm whose entries in the eigenvector basis of the first one are independent with a variance profile.
More explicitly, let $D_n$ be an $n\ti n$ Hermitian matrix, that, up to a change of basis, we suppose diagonal\footnote{If the perturbing matrix belongs to the GOE or  GUE, then its law is invariant under this change of basis, hence our results in fact apply to any self-adjoint matrix   $D_n$.}. We denote by $\mu_n$ the empirical spectral distribution  of $D_n$. This matrix is additively perturbed by a random Hermitian matrix   $\eps_n X_n$ whose entries are chosen at random independently and scaled so that the operator norm of $X_n$ has order one. We are interested in the empirical spectral distribution $\mu_n^\eps$ of  $$ D_n^\eps\;:=\;D_n+\eps_n X_n$$ in the regime where the matrix size $n$ tends to infinity and $\eps_n$ tends to $0$. We shall prove that, depending on the order of magnitude of the perturbation, several regimes can appear. We suppose that $\mu_n$ converges to a limiting measure $\rho(\lam)\ud \lam$ and that the 
variance profile of the entries of $ X_n$  has a macroscopic limit $\sid $ on the diagonal and $\si$ elsewhere.
We then  prove that there is a deterministic function $F$   and a Gaussian random linear form $\ud Z$ on the space of $\Cc^6$ functions on $\R$, both depending only on the limit parameters of the model $\rho,\si$ and $\si_{\op{d}}$   \st if one defines  the distribution   $\ud F:  \phi \longmapsto-\int\phi'(s)F(s)\ud s$, then,  for large $n$:
\begin{align}
\la{densite219161intro}\hspace{4.5cm}
\mu_n^\eps & \;\approx\; \mu_n+\f{\eps_n}{n}\ud Z  &\text{ if }& \eps_n\ll n^{-1}\\
\la{densite2191611intro} \qquad\mu_n^\eps & \; \approx \; \mu_n+\f{\eps_n}{n}\lf(c\ud F + \ud Z \ri) &\text{ if }& \eps_n\sim \f{c}{n}\\
\la{densite2191612intro}\qquad \mu_n^\eps & \; \approx \; \mu_n+\eps_n^2\ud F &\text{ if }& n^{-1}\ll \eps_n\ll 1
\end{align}
and if, moreover, $ n^{-1}\ll \eps_n\ll n^{-1/3},$  then convergence \eqref{densite2191612intro} can be  refined as follows:
\begin{align}
\la{densite2191613intro}\; \mu_n^\eps&\;\approx\; \mu_n+\eps_n^2\ud F+ \f{\eps_n}{n} \ud Z .
\end{align}
In   Section \re{sec:ex}  
several figures show a very good matching of random simulations with  these theoretical results.  The definitions of the function $F$ and of the process $Z$ are given below in \eqre{IntroF} and \eqre{IntroX_phi1}. In many cases, the linear form $\ud F$ can be interpreted as the integration with respect to the signed measure $F'(x)\ud x$. The function $F$ is related to free probability theory, as explained in Section \re{sec:relation_to_free_proba} below, whereas the linear form $\ud Z$ is related to the so-called one-dimensional Gaussian free field defined, for instance, at  \cite[Sect. 4.2]{GFF}. If the variance profile of $X_n$ is constant, then it is precisely the Laplacian of the Gaussian free field, defined in the sense of distributions.

The transition at $\eps_n\sim n^{-1}$ is the well-known transition, in quantum mechanics, where  the \emph{perturbative regime} ends. Indeed, one can distinguish the two following regimes:\bgt\ite The regime $\eps_n\ll n^{-1}$,   called      the \emph{perturbative regime} (see \cite{Fernandez}): the size of the perturbation (\ie its   operator norm)  is  much smaller than the typical spacing between two consecutive eigenvalues (level spacing), which is of order $n^{-1}$ in our setting.\\
\ite The regime  $n^{-1}\ll \eps_n\ll 1$, sometimes called the  \emph{semi-perturbative regime}, where the size of the perturbation  is not small compared to the level spacing. This regime concerns many applications \cite{LedoitPeche,AllezB1} in the context of covariance matrices and applications to finance.\ent 
A surprising fact discovered during this study is that the semi-perturbative regime $n^{-1}\ll \eps_n\ll 1$ decomposes into infinitely many sub-regimes. In the case $n^{-1}\ll \eps_n\ll n^{-1/3}$, the expansion of $\mu_n^\eps-\mu_n$ contains a single deterministic term before the   random term  $\f{\eps_n}{n}\ud Z$. In the case $n^{-1/3}\ll \eps_n\ll n^{-1/5}$, the expansion of $\mu_n^\eps-\mu_n$ contains two of them. More generally, for all positive integer $p$, when $n^{-1/(2p-1)}\ll \eps_n\ll n^{-1/(2p+1)}$, the expansion contains $p$ of them. For computational complexity reasons, the only case we state explicitly is the first one. We refer the reader to Section \ref{subsec:extensions} for a discussion around this point.

 In the   papers  \cite{WW95, AllezB1, AllezB2, AllezB3, AllezB4},     Wilkinson, Walker, Allez, Bouchaud  \emph{et al}   have investigated some     problems related to this one. Some of these works were motivated by     the   estimation of a  matrix out of the observation of its noisy version. Our paper differs from these ones mainly by the facts that firstly, we are interested in  the perturbations of the  \emph{global} empirical distribution of the eigenvalues and not of a single one, and secondly, we push our expansion up to the random term, which  does not appear in these papers. Besides, the noises they consider have constant variance profiles (either a Wigner-Dyson noise in the four first cited papers or a rotationally invariant noise in the fifth
one).  
 The transition at $\eps_n\sim n^{-1}$ between the perturbative and the semi-perturbative regimes   is already present in  these texts. They  also consider the transition between the perturbative regime $\eps_n\ll 1$ and the \emph{non perturbative} regime $\eps_n\asymp  1$.   As explained above,  we exhibit the existence of an infinity of sub-regimes in this transition and focus on $\eps_n\ll 1$ for the first order of the expansion and to $\eps_n\ll n^{-1/3}$ for the second (and last) order. The study of other sub-regimes is    postponed  to forthcoming papers.

The paper is organized as follows. Results, examples and comments are given in Sections    \re{sec:MR} to \re{sec:relation_to_free_proba}, while    the rest of the paper, including an appendix, is devoted to the proofs, except for Section \ref{subsec:extensions}, where we discuss the sub-regimes mentioned above. 

{\bf Notations.} For $a_n,b_n$ some real sequences, $a_n\ll b_n$ (resp. $a_n\sim b_n$)  means that $a_n/b_n$ tends to $0$ (resp. to $1$). Also, $\cvP$ and $\cvd$ stand respectively for convergence in \pro and convergence in distribution for all finite marginals.

\section{Main result}\la{sec:MR}

\subsection{Definition of the model and assumptions}
For all positive integer $n$, we consider   a     real diagonal matrix $D_n =\diag(\lam_n(1), \ld, \lam_n(n))$,  as well as a Hermitian random matrix $$X_n=\ff{\sqrt{n}}[x^n_{i,j}]_{1\le i,j\le n}$$  and a positive number $\eps_n$. The normalizing factor  $ n^{-1/2}$ and our hypotheses below ensure that the operator norm of $X_n$ is of order one.
We then define, for all $n$, $$D_n^\eps:=D_n+ \eps_n X_n.$$ 
We now introduce the probability measures $\mu_n$ and $\mu_n^\eps$ as the respective uniform distributions on the eigenvalues (with multiplicity) of  $D_n$ and $D_n^\eps$. 
 Our aim is to give a perturbative expansion of  $\mu_n^\eps$  around   $\mu_n$.
 
We make the following hypotheses:
 \begin{enumerate}\item[(a)]  the entries $x^n_{i,j}$ of $\sqrt{n}X_n$ are independent (up to symmetry) random variables, centered, with variance denoted  by $\si_n^2(i,j)$, such that $\E   |x^n_{i,j}|^{8} $ is bounded uniformly on $n,i,j$,\\
  \item[(b)] there are $f,\sid,\si$  real functions defined respectively on $[0,1]$, $[0,1]$  and $[0,1]^2$ such that, for each $x\in [0,1]$,   $$\lam_n(\lfl nx\rfl)\ninf f(x)\quad\trm{ and }\quad \si_n^2(\lfl nx\rfl,\lfl nx\rfl)\ninf \sid(x)^2$$ and for each $x\ne y\in [0,1]$,
  $$\si_n^2(\lfl nx\rfl,\lfl ny\rfl)\ninf \si^2(x,y).$$
We  make the following hypothesis about the rate of convergence:
$$\eta_n\;:=\;\max\{n\eps_n,1\}\ti   \sup_{1\le i\ne j\le n}( |\si_n^2(i,j)-\si^2(i/n,j/n)|+|\lam_n(i)-f(i/n)|)\ninfe 0.$$  
  \end{enumerate}
  Let us now make some   assumptions on the limiting functions  $\si$ and $f$: \begin{enumerate}\item[(c)] the function $f$ is bounded and the push-forward   of the uniform measure on $[0,1]$ by the function $f$ has a  density $\rho$ with respect to the Lebesgue measure on $\R$ and a   compact support denoted by $\Sc$,\\
 \item[(d)] the variance of the entries of $X_n$ essentially depends on the eigenspaces of $D_n$, namely, there exists a symmetric    function $\tau(\sd,\sd )$ on $\R^2$ \st for all $x\ne y$, $\si^2(x,y)=\tau(f(x),f(y))$,
   \\
 \item[(e)] the following regularity property holds: there exist $\eta_0>0, \al>0$ and $C<\infty$  \st   for almost all $s\in \R$, for all $t\in [s-\eta_0, s+\eta_0]$,\quad $|\tau(s,t)\rho(t)-\tau(s,s)\rho(s)|\le C|t-s|^\al$.
\end{enumerate}
We add a last assumption which strengthens assumption (c) and makes it possible to include the case where the set of eigenvalues of $D_n$ contains some outliers:
\begin{enumerate}
\item[(f)] there is a real compact set $\tS$   \st $$\max_{1\le i\le n}\op{dist}(\lam_n(i), \tS)\ninf 0.$$
\end{enumerate}

\beg{rmq}[About the hypothesis that $D_n$ is diagonal] (i) If the perturbing matrix $X_n$ belongs to the GOE (resp. to the  GUE), then its law is invariant under   conjugation by any orthogonal (resp. unitary) matrix. It follows that in this case,   our results   apply  to any real symmetric (resp. Hermitian) matrix  $D_n$ with eigenvalues $\lam_n(i)$ satisfying the above hypotheses.
\ite (ii) As explained after Proposition \re{prop:relation_toSch96} below, we conjecture that when the variance profile of $X_n$ is constant, for $\eps_n\gg n^{-1}$, we do not need the hypothesis that $D_n$ is diagonal neither. 
However, if the perturbing matrix does not have a constant variance profile,   then for a non-diagonal $D_n$ and $\eps\gg n^{-1}$, the spectrum of $D_n^\eps$ should depend heavily on the relation between the eigenvectors of $D_n$ and the variance profile, which implies that our results should not remain true. 
\ite (iii) At last,  it is easy to see that the random process $(Z_\phi)$ introduced at \eqre{IntroX_phi1} satisfies, for any test function $\phi$, 
$$\ff{\eps_n}\sum_{i=1}^n\lf(\phi(\lam_n(i)+\f{\eps_n}{\sqrt n}x_{ii})-\phi(\lam_n(i))\ri)\ninfd Z_\phi.$$
 Thus, regardless to the variance profile,  the convergence of \eqre{1731617h12} rewrites, informally, 
\be\la{229171} \mu_n^\eps=\ff n\sum_{i=1}^n\delta_{\lam_n(i)+(\eps_n/\sqrt n)x_{ii}}+o(\eps_n/n).\ee A so simple expression, up to a $o(\eps_n/n)$ error, of the empirical spectral distribution of $D_n^\eps$, with some independent translations  $\f{\eps}{\sqrt n}x_{ii}$,  should not remain true without the hypothesis that $D_n$ is diagonal or that 
the distribution of $X_n$ is invariant under conjugation. 

\en{rmq}

\subsection{Main result}

Recall that the {\it Hilbert transform}, denoted by $H[u]$,  of a   function $u$, is the function $$H[u](s):=\operatorname{p.v.}\int_{t\in\R}\f{u(t)}{s-t}\ud t$$
 and define the function    \be\la{IntroF}F(s) = -\rho(s) H[\tau(s,\cdot)\rho(\cdot)]
(s).\ee 
 Note that, by assumptions (c) and (e), $F$ is well defined and supported by $\Sc$.  Besides, for any $\phi$ supported on an interval where $F$ is $\Cc^1$, $$-\int\phi'(s)F(s)\ud s=\int \phi(s)\ud F(s),$$ where $\ud F(s)$ denotes the measure $F'(s)\ud s$.

We also introduce the centered Gaussian field, $(Z_\phi)_{\phi\in \Cc^6}$, indexed by the set of $\Cc^6$ complex functions on $\R$, with covariance defined by
 \be\la{IntroX_phi1}\E Z_\phi  Z_\psi \;=\;\int_0^1\sid(t)^2\phi'(f(t))  \psi'(f(t))\ud t\qquad \trm{ and }\qquad \ovl{Z_\psi}\;=\;  Z_{\ovl{\psi}}
.\ee
Note that the process $(Z_\phi)_{\phi\in\mathcal{C}^6}$ can be represented, for $(B_t)$ is the standard one-dimensional Brownian motion,  as
$$Z_\phi = \int_0^1 \sid(t) \phi'(f(t)) \ud B_t.$$

\beg{Th}\label{ConvergenceTheorem00} For all compactly supported $\Cc^6$ function $\phi$ on $\R$, the following convergences hold:

$\bullet$ {\bf Perturbative regime:}  if $\eps_n\ll n^{-1}$, then, \be\la{1731617h12} n\eps_n^{-1}(\mu_n^\eps-\mu_n)(\phi) \ninfd Z_\phi.\ee
$\bullet$ {\bf Critical regime:}  if $\eps_n\sim c/n$, with $c$ constant, then, \be\la{1731617h13} n\eps_n^{-1}(\mu_n^\eps-\mu_n)(\phi) \ninfd -c\int \phi'(s)F(s)\ud s+Z_\phi.\ee
$\bullet$  {\bf Semi-perturbative regime:} if $  n^{-1}\ll \eps_n\ll 1$, then, 
    \be\la{1731617h14} \eps_n^{-2}(\mu_n^\eps-\mu_n)(\phi)\ninfP  -\int \phi'(s)F(s)\ud s,\ee and if, moreover,  
    $  n^{-1}\ll \eps_n\ll n^{-1/3}$, then,
   \be\la{1731617h15} n\eps_n^{-1}\lf((\mu_n^\eps-\mu_n)(\phi) + \eps_n^2\int \phi'(s)F(s)\ud s\ri) \ninfd Z_\phi.\ee
\en{Th}

\begin{rmq}[Sub-regimes for $n^{-1/3}\ll\eps_n\ll 1$]\la{caselargereps}
In the semi-perturbative regime, the reason why we provide an expansion up to a random term, only for   $\eps_n\ll n^{-1/3}$, is that the study of the regime $n^{-1/3}\ll\eps_n\ll1$ up to such a precision, requires further terms in the expansion of the resolvent of $D_n^\eps$ that make appear, beside $\ud F$, additional determistic terms of smaller order, which are much larger than the probabilistic term containing $Z_\phi$.  The computation becomes rather intricate without any clear recursive formula.
As we will see in Section \ref{subsec:extensions},   there are infinitely many regimes.    Precisely, for any positive integer $p$, when $n^{-1/(2p-1)}\ll \eps_n\ll n^{-1/(2p+1)}$, there are $p$ deterministic terms in the expansion before the term in $Z_\phi$. 
\end{rmq}

\begin{rmq}[Local law]\la{rmk:DL22916} 
The  approximation $$\ds \mu_n^\eps(I)\;\approx \; \mu_n (I)+ \eps_n^2\int_I\ud F$$ of \eqre{1731617h14}   should stay true even for intervals $I$ with size tending to $0$ as the dimension $n$ grows, as long as the size of $I$ stays much larger than the right-hand side term of \eqre{71171}, as can be seen from Proposition \re{localLemma}. 
\end{rmq}

\beg{rmq}\la{hypDiscNath}The second part of Hypothesis (b), concerning the speed of convergence of the profile of the spectrum of $D_n$ as well as of the variance of its perturbation, is needed 
in order to express the expansion of $\mu_n^\eps-\mu_n$ in terms of 
 limit parameters of the model $  \si$ and $\rho$.
We can   remove this hypothesis and get  analogous expansions where the terms $\ud F$ and $\ud Z$ are replaced by  their discrete counterparts   $\ud F_n$ and $\ud Z_n$, defined thanks to the ``finite $n$" empirical versions  of the limit parameters   $\si$ and $\rho$. 
\en{rmq}

\section{Examples}\la{sec:ex}
\subsection{Uniform measure perturbation by a band matrix}\label{example1}
Here, we consider the case where $f(x)=x$, $\si_{\op{d}}(x)\equiv  m$  and $\si(x,y)=\one_{|y-x|\le \ell}$, for   some   constants $m\ge 0$ and $\ell\in [0,1]$ (the relative width of the band).
In this case, $\tau(\sd,\sd)=\si(\sd,\sd)^2$, hence   \be\la{10mai16}F(s)\;=\;\one_{(0,1)}(s) \op{p.v.}\int_{t}\f{\tau(s,t)}{s-t}\ud t
\;=\;- \one_{(0,1)}(s)\log \f{\ell \wedge(1-s)}{\ell\wedge s} \ee
and $(Z_\phi)_{\phi\in \Cc^6}$ is the  centered complex Gaussian process   with covariance defined by
$$
\mathbb{E}Z_\phi \ovl{Z_\psi} \; =\; m^2 \int_0^1 \phi'(t) \ \ovl{\psi'(t)}  \ \ud t
\qquad \trm{ and }\qquad \ovl{Z_\psi}\;=\;  Z_{\ovl{\psi}}.
$$

\begin{figure}[h!]
\centering
\subfigure[$n=10^4$, $\eps_n=n^{-0.4}$,   $\ell=0.2$]{
\includegraphics[width=2.45in]{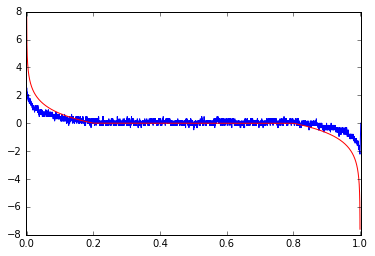}\label{unif_width0.2}}\hspace{1.2in}
\subfigure[$n=10^4$, $\eps_n=n^{-0.4}$,  $\ell=0.8$]{
\includegraphics[width=2.45in]{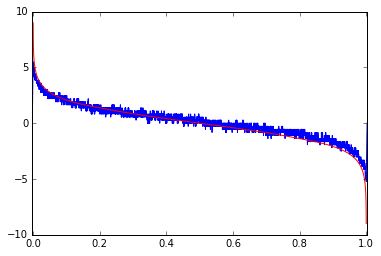}
\label{unif_width0.8}}
\caption{{\bf Deforming the uniform distribution by a   band matrix.} Cumulative   distribution function  of $\eps_n^{-2}(\mu_n^\eps-\mu_n)$  
 (in blue) and   function $F(\sd)$ of \eqre{10mai16} (in red). The non smoothness of the blue curves results of the noise term $Z_\phi$ in Theorem  \re{ConvergenceTheorem00}.
 Each graphic is realized thanks to one single matrix (no averaging) perturbed by a real Gaussian band matrix.}
\label{fig:unif_band}
\end{figure}

Theorem \re{ConvergenceTheorem00}   is then  illustrated by Figure \ref{fig:unif_band}, where we ploted the cumulative distribution functions. 

\subsection{Triangular pulse perturbation by a Wigner matrix}
\begin{figure}[h!]
\begin{center}
\includegraphics[scale=.4]{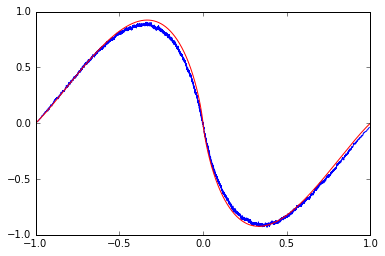}
\includegraphics[scale=.4]{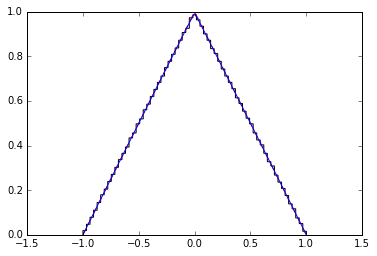}
\includegraphics[scale=.4]{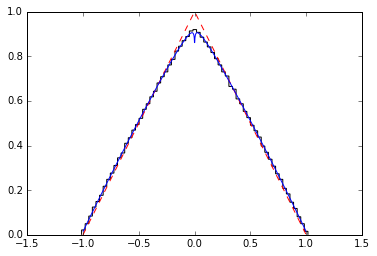}
\includegraphics[scale=.4]{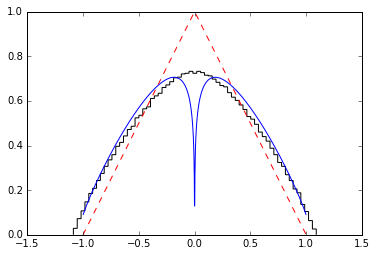}
\caption{{\bf Triangular pulse perturbation by a Wigner matrix: density and cumulative distribution   function.} \emph{Top left:}  Cumulative   distribution function  of $\eps_n^{-2}(\mu_n^\eps-\mu_n)$(in blue) and   function $F(\sd)$ of \eqre{FTriPulse22916} (in red).
\emph{Top right and bottom:}    Density $\rho$  (red dashed line), histogram of the eigenvalues of $D_n^\eps$ (in black)  and     theoretical density $\rho+\eps_n^2F'(s)$ of the eigenvalues of $D_n^\eps$ as predicted by Theorem \re{ConvergenceTheorem00}  (in blue).
Here,    $n= 10^4$ and  $\eps_n=n^{-\alpha}$, with $\alpha=0.25$  (up left),  $\alpha=0.4$ (up right), $0.25$ (bottom left) and $0.1$   (bottom  right).}\label{fig:tri_pluse}
\end{center}
\end{figure}

Here, we consider the case where $\rho(x)=(1-|x|)\one_{[-1,1]}(x)$, $\si_d \equiv m$, for some real constant $m$, and $\si\equiv 1$ (what follows can be adapted to the case $\si(x,y)=\one_{|y-x|\le \ell}$, with a bit longer formulas).
In this case, thanks to the formula (9.6)  of $H[\rho(\sd)]$ given p. 509 of \cite{kingvol2}, we get \be\la{FTriPulse22916}F(s)\;=\;(1-|s|)\one_{[-1,1]}(s)\lf\{(1-s)\log(1-s)-(1+s)\log(1+s)+2s\log |s|\ri\}.\ee
and the  covariance of  $(Z_\phi)_{\phi\in \Cc^6}$ is given by  
$$
\mathbb{E}Z_\phi \ovl{Z_\psi} \; =\; m^2 \int_{-1}^1 (1-|t|) \  \phi'(t) \ \ovl{\psi'(t)} \  \ud t
\qquad \trm{ and }\qquad \ovl{Z_\psi}\;=\;  Z_{\ovl{\psi}}.
$$
Theorem \re{ConvergenceTheorem00}   is then illustrated by Figure \ref{fig:tri_pluse} in the case where $\eps_n\gg n^{-1/2}$. In Figure \ref{fig:tri_pluse}, we implicitly use some test functions of the type $\phi(x)=\one_{x\in I}$ for some intervals $I$. These functions are not $\Cc^6$, and one can easily see that for $\eps_n\ll n^{-1/2}$,  Theorem \re{ConvergenceTheorem00} cannot work for such functions. However, considering imaginary parts of  Stietljes transforms, \ie test functions $$\ds \phi(x)=\ff{\pi}\f{\eta}{(x-E)^2+\eta^2}\qquad (E\in \R, \; \eta>0)$$ gives a perfect matching between the predictions from  Theorem  \re{ConvergenceTheorem00} and numerical simulations, also for $\eps_n\ll n^{-1/2}$ (see Figure \re{fig:triangular_pulse}, where we use Proposition \re{StieltjesTheorem} and  
\eqre{032016112} to compute the theoretical limit). 
\begin{figure}[h!]
\centering
\includegraphics[width=1.8in]{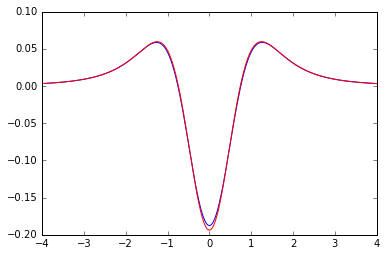}
\includegraphics[width=1.8in]{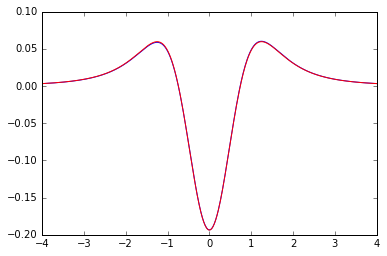}
\includegraphics[width=1.8in]{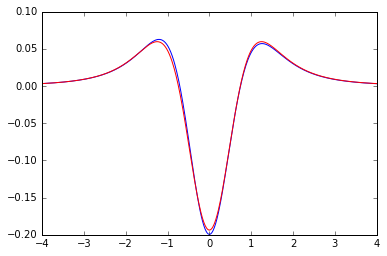} 
\caption{{\bf Triangular pulse perturbation by a Wigner matrix:   Stieltjes transform.} Imaginary part of the Stieltjes transform of $\eps_n^{-2}(\mu_n^{\eps} - \mu_n)$ (in blue) and of the measure $\ud F$ (in red) at $z=E+\ii$ as a function of the real part $E$  for different values of $\eps_n$. Here, $n= 10^4$ and  $\eps_n=n^{-\alpha}$, with $\alpha=0.2$, $0.5$ and $0.8$ (from left to right).}
\label{fig:triangular_pulse}
\end{figure}

\subsection{Parabolic    pulse perturbation by a Wigner matrix}
\begin{figure}[h!]
\begin{center}
\includegraphics[scale=.4]{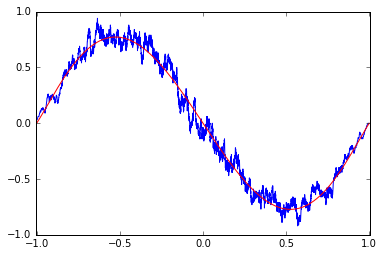}
\includegraphics[scale=.4]{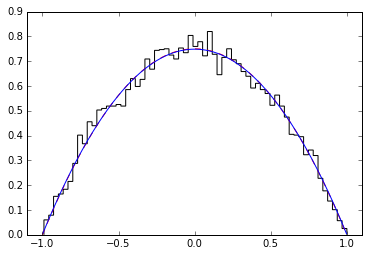}
\includegraphics[scale=.4]{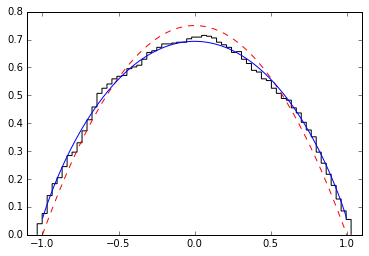}
\includegraphics[scale=.4]{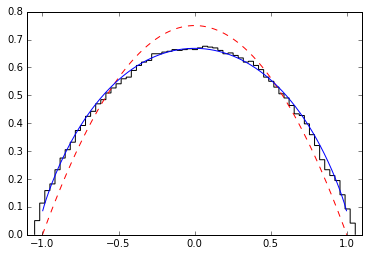}
\caption{{\bf  Parabolic pulse perturbation  by  a Wigner matrix.}     \emph{Top left:}  Cumulative   distribution function  of $\eps_n^{-2}(\mu_n^\eps-\mu_n)$(in blue) and   function $F(\sd)$ of \eqre{quadraticF22916} (in red).
\emph{Top right and bottom:}    Density $\rho$  (red dashed line), histogram of the eigenvalues of $D_n^\eps$ (in black)  and     theoretical density $\rho+\eps_n^2F'(s)$ of the eigenvalues of $D_n^\eps$ as predicted by Theorem \re{ConvergenceTheorem00}  (in blue).
Here,   $n= 10^4$ and  $\eps_n=n^{-\alpha}$, with $\alpha=0.25$  (up left),  $\alpha=0.4$ (up right), $0.2$ (bottom left) and $0.18$   (bottom  right).}\label{fig:para_pluse}
\end{center}
\end{figure}

Here, we consider the case where $\rho(x)=\f{3}{4}(1-x^2)\one_{[-1,1]}(x)$, $\si_d \equiv m$, for some real constant $m$, and $\si\equiv 1$ (again, this can be adapted to the case $\si(x,y)=\one_{|y-x|\le \ell}$). Theorem \re{ConvergenceTheorem00}   is then illustrated by Figure \ref{fig:para_pluse}.
In this case, thanks to the formula (9.10)  of $H[\rho(\sd)]$ given p. 509 of \cite{kingvol2}, we get \be\la{quadraticF22916}F(s)\;=\; -\f{9}{16}(1-s^2)\one_{[-1,1]}(s)\lf\{ 2s - (1-s^2) \ln \lf| \f{s-1}{s+1} \ri| \ri\}\ee
and the  covariance of  $(Z_\phi)_{\phi\in \Cc^6}$ is given by  
$$
\mathbb{E}Z_\phi \ovl{Z_\psi} \; =\; \f{3m^2}{4} \int_{-1}^1 (1-t^2) \ \phi'(t) \ \ovl{\psi'(t)} \  \ud t\qquad \trm{ and }\qquad \ovl{Z_\psi}\;=\;  Z_{\ovl{\psi}}.$$

 \section{Relation to free \pro theory}\la{sec:relation_to_free_proba}
 
 Let us now explain how this work is related to free probability theory.  If, instead of letting    $\eps_n$ tend to zero, one considers the model $$D_n^t\;:=\;D_n+ \sqrt{t} X_n$$ for a fixed $t>0$, then, by  \cite{CasatiGirko1,CasatiGirko2,Sh96,AZ05}, the empirical eigenvalue distribution of $D_n^t$ has a limit as $n\to\infty$, that we shall denote here by $\mu_t$. The law $\mu_t$ can be interpreted  as the law of the sum of two elements in a non-commutative \pro space which are free with an amalgamation over a certain sub-algebra (see \cite{Sh96} for more details).  The following proposition relates the function $F$ from \eqre{IntroF} to the  first order expansion of $\mu_t$ around $t=0$. 
 \beg{propo}\la{prop:relation_toSch96}
For any $z\in \C\bck\R$, we have  $$ \f{\partial}{\partial t}_{|t=0}  \int  \f{\ud\mu_t(\lam)}{z-\lam}\;=\;-\int\f{F(\lam)}{(z-\lam)^2}\ud \lam\;=\;-\int F(\lam)  \f{\partial}{\partial \lam}\lf(\ff{z-\lam} \ri)\ud\lam.$$
  \en{propo}
 
This is related to the fact  that in  Equations \eqre{densite219161intro}--\eqre{densite2191613intro}, for  $\eps_n$ large enough, the term $\eps_n^2\ud F$ is the leading term. 

In the particular case where $X_n$ is a Wigner matrix, $\mu_t$ is the free convolution of the measure $\rho(\lam)\ud \lam$ with a semicircle distribution and admits a density  $\rho_t$, by \cite[Cor. 2]{Biane2}. 
Then,     
Theorem \ref{ConvergenceTheorem00} 
makes it possible to formally  recover  the \emph{free Fokker-Planck equation with null potential}:
\bes\label{burgersreel}\begin{cases}\pt \rho_t(s)+\ps\{\rho_t(s)H[\rho_t ](s)\}=0,\\
\rho_0(s)=\rho(s),\end{cases}\ees where $H[\rho_t ]$  denotes the Hilbert transform of $\rho_t$. This equation is 
    also called  \emph{McKean-Vlasov} (or \emph{Fokker-Planck}) \emph{equation  with logarithmic interaction} (see \cite{Biane3,Biane4,Biane5}).  
    
    
    Note also that when $X_n$ is a Wigner matrix, the hypothesis that $D_n$ is diagonal is not required to have the convergence of the empirical eigenvalue distribution of $D_n^t$ to $\mu_t$ as $n\to\infty$. This suggests that, even for non diagonal $D_n$,  the convergence of  \eqre{1731617h14} still holds when $X_n$ is a Wigner matrix.

  {\bf Proof of Proposition \re{prop:relation_toSch96}.}
   By  \cite[Th. 4.3]{Sh96}, we have  \be\la{8jan20169h} \int  \f{\ud\mu_t(\lam)}{z-\lam}\;=\;\int_{x=0}^1C_t(x,z)\ud x,\ee where $C_t(x,z)$  is bounded by $|\Im z|^{-1}$ and satisfies the fixed-point equation \bes\label{eqC_epsilon}C_t(x,z)=\ff{z-f(x)-t\int_{y=0}^1 \si^2(x,y)C_t(y,z)\ud y}.\ees
Hence   as $t\to 0$,   
$C_t(x, z)\lto \ff{z-f(x)}$ uniformly in $x$.
Thus \beq C_t(x,z)-\ff{z-f(x)}&=&\f{t\int_{y=0}^1\si^2(x,y)C_t(y,z)\ud y}{(z-f(x)-t\int_{y=0}^1\si^2(x,y)C_t(y,z)\ud y)(z-f(x))}\\ 
&=&
t\ff{(z-f(x))^2} \int_{y=0}^1\si^2(x,y)C_t(y,z)\ud y+o(t)\\
&=&
t\ff{(z-f(x))^2} \int_{y=0}^1\f{\si^2(x,y)}{z-f(y)}\ud y+o(t)
\eeq
where each $o(t)$   is   uniform in $x\in [0,1]$. Then, by \eqre{8jan20169h}, we deduce that  $$ \f{\partial}{\partial t}_{|t=0}  \int  \f{\ud\mu_t(\lam)}{z-\lam}= \int_{(x,y)\in [0,1]^2} \f{\si^2(x,y)}{(z-f(x))^2(z-f(y))}\ud x\ud y .$$
The right-hand side term of the previous equation is precisely the number $B(z)$ introduced at \eqre{032016112} below. Then, one concludes using Proposition \re{StieltjesTheorem} from Section \re{sec:5123}.
 \hfill$\square$
 
 \section{Strategy of the proof}\la{sec:strategy}

We shall first prove the convergence results of Theorem \ref{ConvergenceTheorem00} for test functions $\phi$ of the form  $\varphi_z(x):=\ff{z-x}$. This is done in Section \ref{sec:stieltjes} by writing an expansion of the resolvent of $D_n^\eps$.

Once we have proved that the convergences hold for the resolvent of $D_n^\eps$, we can extend them to the larger class of compactly supported $\mathcal{C}^6$ functions on $\R$.

In Section \ref{sec:probability}, we use the Helffer-Sj\"ostrand formula to extend the convergence in probability in the semi-pertubative regime \eqref{1731617h14} to the case of compactly supported $\mathcal{C}^6$ functions on $\R$.

In Section \ref{sec:clt}, the convergences in distribution \eqref{1731617h12}, \eqref{1731617h13} and \eqref{1731617h15} are proved in two steps. The overall strategy is to apply an extension lemma of Shcherbina and Tirozzi  which states that a CLT  that applies to a sequence of \textit{centered} random linear forms on some   space   can be extended, by density, to a larger   space, as long as the variance of the image of these random linear forms by a function $\phi$ of the larger space is uniformly bounded by the norm of $\phi$. Therefore, our task is twofold.
We need first to prove that the sequences of variables involved in the convergences \eqref{1731617h12}, \eqref{1731617h13} and \eqref{1731617h15} can be replaced by their \textit{centered} counterparts $n\eps_n^{-1}(\mu_n^\eps(\phi)-\E[\mu_n^\eps(\phi)])$ (\ie they differ by $o(1)$). In a second step, we dominate the variance of these latter variables, in order to apply the extension lemma which is precisely stated in the appendix as Lemma \ref{extensionLemma}.

 \section{Stieltjes transforms convergence}\la{sec:stieltjes}

As announced in the previous section, we start with the proof of Theorem \ref{ConvergenceTheorem00} in the special case of test functions of the type $\varphi_z:= \ff{z-x}$. We decompose it into two propositions. Their statement and proof are the purpose of the three following subsections. The two last subsections \ref{subsec:local} and \ref{subsec:extensions} are devoted  respectively to a local type convergence result and to a discussion about the possibility of an extension of the expansion result to a wider range of rate of convergence of $\eps_n$, namely beyond 
$n^ {-1/3}$. 
 
 \subsection{Two statements}\la{sec:5123}

 Let denote, for 
 $z\in\C\bck\R$, 
 \be\la{03201611}Z(z):=Z_{\varphi_z}\quad \trm{ for }\quad\varphi_z(x):=\ff{z-x}\ee where $(Z_\phi)_{\phi\in \Cc^6}$ is the Gaussian field with covariance defined by \eqre{IntroX_phi1}.
 We also introduce, for  $z\in\C\bck\R$,    \be\la{032016112}B(z):=\int_{(s,t)\in [0,1]^2}\f{\si^2(s,t)}{(z-f(s))^2(z-f(t))}\ud s\ud t\ee
 and   \be\la{032016113}\DG_n(z)\;:=\; (\mu_n^\eps-\mu_n)(\varphi_z)\;=\;\ff{n}\Tr\ff{z-D_n^\eps}-\ff{n}\Tr\ff{z-D_n}.\ee

\beg{propo}\label{ConvergenceTheorem} Under Hypotheses (a), (b), (f),
\bgt\ite if $\eps_n\ll n^{-1}$, then for all $z\in\C\bck\R$,
\be\la{18mars161} n\eps_n^{-1}\DG_n(z) \ninfd Z(z)\,\ee
   \ite if $\eps_n\sim c/n$, with $c$ constant, then for all $z\in\C\bck\R$
\be\la{18mars162} n\eps_n^{-1}\DG_n(z) \ninfd cB(z)+Z(z) \, ,\ee
   \ite if $  n^{-1}\ll \eps_n\ll n^{-1/3}$,   then for all $z\in\C\bck\R$
    \be\la{18mars164} n\eps_n^{-1} \lf(\DG_n(z)-\eps_n^2B(z)\ri) \ninfd  Z(z)\,.\ee   
    \ite if $  n^{-1}\ll \eps_n\ll 1$, then for all $z\in\C\bck\R$,
    \be\la{18mars166}\eps_n^{-2}\DG_n(z) - B(z) \ninfP 0 \,.\ee
\ent
\en{propo}

{\noindent{\bf{Remark.}} Note that \eqre{18mars162} is merely an extension of \eqre{18mars164} in the critical regime.}

The following statement expresses $B(z)$ as the image of a $\varphi_z$ by a linear form. So, in the expansion of the previous proposition, both quantities $Z(z)$ and $B(z)$ depend linearly on $\varphi_z$. Note that as $F$ vanishes at $\pm\infty$, Proposition \re{StieltjesTheorem} does not contradicts the fact that as $|z|$ gets large, $B(z)=O(|z|^{-3})$.

\beg{propo}\la{StieltjesTheorem}Under Hypotheses (c), (d), (e), for any $z\in \C\bck\Sc$, for $F$ defined by \eqre{IntroF},  
$$B(z)\;=\;-\int\f{F(s)}{(z-s)^2}\ud s\;=\; -\int \varphi_z'(s)F(s)\ud s.$$
\en{propo} 

\subsection{Proof of Proposition \ref{ConvergenceTheorem}}

 The proof is based on a perturbative expansion of the resolvent $\ff{n} \Tr\ff{z-D_n^\eps}$. To make notations lighter, we shall sometimes  suppress the subscripts and superscripts $n$, so that $D_n^\eps$, $D_n$, $X_n$ and $x_{i,j}^n$ will be respectively denoted by $D^\eps$, $D$, $X$ and $x_{i,j}$.
  Let us fix $z\in \C\bck\tS$.   
    We can deduce from the expansion of the resolvent of $D^\eps$:
 \bes\la{131015}\DG_n(z) = A_n(z)+B_n(z)+C_n(z)+R_n^\eps(z), \ees
  with \beq\la{formuleA} A_n(z)&:=&  \f{\eps_n}{n} \Tr\ff{z-D}X\ff{z-D}=    \f{\eps_n}{n}  \ff{\sqrt{n}}\sum_{i=1}^n\f{x_{i,i}}{(z-\lam_n(i))^2}\\
\la{formuleB}  B_n(z)&:=& \f{\eps_n^2}{n}\Tr \ff{z-D}X\ff{z-D}X\ff{z-D}=\f{\eps_n^2}{n^2}\sum_{i,j}\f{|x_{i,j}|^2}{(z-\lam_n(i))^2(z-\lam_n(j))} \\
\la{formuleC}  C_n(z)&:=&\f{\eps_n^3}{n}\Tr \ff{z-D} X \ff{z-D}X\ff{z-D}X \ff{z-D} \nonumber \\ &=& \f{\eps_n^3}{n^{5/2}} \sum_{i,j,k=1}^n \f{x_{i,j} \ x_{j,k} \ x_{k,i}}{(z-\lambda_n(i))^2 \ (z-\lambda_n(j)) \ (z-\lambda_n(k))} 
\\
\la{formuleR}  R_n^\eps(z)&:=&\f{\eps_n^4}{n}\Tr \ff{z-D} X \ff{z-D} X \ff{z-D}X\ff{z-D}X \ff{z-D^\eps}.
\eeq 

The purpose of the four following claims is to describe the asymptotic behavior of each of these four terms.

\beg{Claim}\la{Claim1} The finite dimension marginals of the centered  process $$(n\eps_n^{-1} A_n(z))_{z\in \C\bck\tS}$$ converge in distribution to those of the    centered Gaussian process $(Z(z))_{z\in \C\bck\tS}$. Besides, there is $C>0$ \st  for any $z\in \C\bck\tS$, \be\la{Eq:VarA(z)}\E[ |n\eps_n^{-1}A_n(z)|^2]\;\le \; \f{C}{\op{dist}(z, \tS)^4}.\ee
\en{Claim}

\bpr  Estimate \eqre{Eq:VarA(z)} follows from 
\bes\la{Eq2:VarA(z)230916}\E[ |A_n(z)|^2] \; = \; \f{\eps_n^{2}}{n^3} \sum_{i=1}^n \f{\mathbb{E}\lf[|x_{i,i}|^2\ri]}{|z-\lambda_n(i)|^4} \; \le \; \f{\eps_n^{2}}{n^3} \sum_{i=1}^n \f{\sigma_n^2(i,i)}{\op{dist}(z,\tS)^4}\ees
and from the existence of a uniform upper bound for $\sigma_n^2(i,i)$ which comes from Hypothesis (a) which stipulates that the $8$-th moments of the entries $x_{i,j}$ are uniformly bounded.

We turn now to the proof of the convergence in distribution of $n\eps_n^{-1}A_n(z)$ which actually does not depend on the sequence $(\eps_n)$. For all $\alpha_1,\beta_1,\dots,\alpha_p,\beta_p \in \C$ and for all $z_1,\dots,z_p \in \C\bck\tS$,
\[\sum_{i=1}^p \alpha_i \left(n\eps_n^{-1}    A_n(z_i) \right) + \beta_i  \overline{\left( n\eps_n^{-1}    A_n(z_i) \right)} = \ff{\sqrt{n}} \sum_{j=1}^n x_{j,j}   \left( \sum_{i=1}^p \xi_n(i,j) \right)
\]
for $\displaystyle{ \xi_n(i,j) = \f{\alpha_i}{(z_i-\lam_n(j))^2} + \f{\beta_i}{(\overline{z_i}-\lam_n(j))^2} }$.  

On one hand, by dominated convergence, the covariance matrix of the 
above two dimensional random vector converges.

On the other hand,   $\E|x_{i,j}|^4$ is uniformly bounded in $i$, $j$ and $n$, by Hypothesis (a). 
Moreover, for $n$ large enough, for all $i,j$,  $$|\xi_n(i,j)|\le 2\max_{1\le i\le p}(|\al_i|+|\beta_i|)\ti (\min_{1\le i\le p} \op{dist}(z_i, \mc{S}))^{-1}.$$
Hence, the conditions of Lindeberg Central Limit Theorem are satisfied and the finite dimension marginals of the process $(n\eps_n^{-1} A_n(z))_{z\in \C\bck\tS}$ converge in distribution to those of the    centered Gaussian process $(Z_z)_{z\in \C\bck\tS}$ defined by its covariance structure 
\beq \E\left(Z(z) \ovl{Z({z'}})\right) &=&\lim_{n\to \infty}\mathbb{E}\left[ \left(n\eps_n^{-1}  A_n(z)\right) . \left(n\eps_n^{-1} \ovl{A_n(z')}\right) \right]\\
&=&\lim_{n\to \infty} \ff{n} \sum_{i,j = 1}^n \f{\mathbb{E}\left[ x_{i,i} \ \ovl{x_{j,j}} \right]}{(z-\lam_n(i))^2 \ (\ovl{z'}-\lam_n(j))^2}\\
&=&  \int_0^1 \f{\sid (t)^2}{(z-f(t))^2 \ (\ovl{z'}-f(t))^2} \ud t \eeq and by the fact that $\overline{Z(z)} = Z(\overline{z})$ which comes from $\overline{A_n(z)} = A_n(\overline{z})$.
\epr

\beg{Claim}\la{Claim2}There is a constant $C$ such that, for $\eta_n$ as in Hypotheses (b),
\bgt\ite if $ \eps_n\ll n^{-1}$, then 
\beq\label{eq:varB(z)1} \E [|n\eps_n^{-1}B_n(z)|^2]\;\le\;  \f{C(n\eps_n)^2}{\op{dist}(z, \tS)^6}+ \f{C\eta_n^2}{\op{dist}(z, \tS)^8},
\eeq
\ite if 
$\eps_n\sim c/n$ or if $ n^{-1} \ll \eps_n\ll 1$, then   
\beq\label{eq:varB(z)2}
\quad \E[|  n\eps_n^{-1}(B_n(z)-\eps_n^2B(z))|^2]\;\le \;\f{C  \eps_n^2}{\op{dist}(z, \tS)^6}+\f{C \eta_n^2 }{\op{dist}(z, \tS)^8}.
\eeq
\ent \end{Claim}

\bpr Remind that,
$$B_n(z)=\f{\eps_n^2}{n^2}\sum_{i,j}\f{|x_{i,j}|^2}{(z-\lam_n(i))^2(z-\lam_n(j))}.$$ 
Introduce the variable $\Bo(z)$ obtained by centering the variable $n\eps_n^{-2} B_n(z)$: $$\Bo(z):=n\eps_n^{-2}(B_n(z)-\E B_n(z))=\ff{n}\sum_{i,j}\f{|x_{i,j}|^2-\si_n^2(i,j)}{(z-\lam_n(i))^2(z-\lam_n(j))}$$ and the defect variable \beqy\la{IntroBodeltan}  \delta_n(z)&:=&\eps_n^{-2}\E B_n(z)-B(z)\nonumber\\ &=&\ff{n^2}\sum_{i,j}\f{\si_n^2(i,j)}{(z-\lam_n(i))^2(z-\lam_n(j))}-\int_{(s,t)\in [0,1]^2}\f{\si^2(s,t)}{(z-f(s))^2(z-f(t))}\ud s\ud t.\nonumber\eeqy
In the two regimes $\eps_n \ll n^{-1}$ and $\eps_n \geq c/n$, we want to dominate the $L^2$ norms respectively of
\[ n\eps_n^{-1}B_n(z)=\eps_n\Bo(z)+n\eps_n(\delta_n(z)+B(z))\quad\text{ and }\quad n\eps_n^{-1}(B_n(z)-\eps_n^2B(z))=\eps_n\Bo + n\eps_n\delta_n(z).
\]
For this purpose, we successively dominate $\Bo$, $\delta_n(z)$ and $B(z)$.

Using the independence of the $x_{i,j}$'s, the fact that they are bounded in $L^4$ and the fact that $z$ stays at a macroscopic distance of the $\lam_n(i)$'s, we can write  for all $z\in\C\bck \tS$
\begin{align}\label{211220161808}
 \E [| \Bo(z)|^2] &=\ff{n^2}\Var\lf(\sum_{  i\le j} \lf(x_{i,j}^2+\one_{i\ne j}\ovl{x_{i,j}}^2\ri)\ff{(z-\lam_n(i))^2(z-\lam_n(j))}\ri)\nonumber\\
&=\ff{n^2}\sum_{  i\le j} \Var\lf(\lf(x_{i,j}^2+\one_{i\ne j}\ovl{x_{i,j}}^2\ri)\ff{(z-\lam_n(i))^2(z-\lam_n(j))}\ri)\nonumber\\
&\le C\op{dist}(z, \tS)^{-6}\, .
\end{align}
Now, the term  $\delta_n(z)$ rewrites
\beq\label{211220161833}  \delta_n(z)
&=&O(n^{-1})\\ \nonumber&&+\int_{(s,t)\in [0,1]^2}\one_{\lfl ns\rfl\ne\lfl nt\rfl}\lf(\f{\si_n^2(\lfl ns\rfl,\lfl nt\rfl)}{(z-\lam_n(\lfl ns\rfl))^2(z-\lam_n(\lfl nt\rfl))}-\f{\si^2(s,t)}{(z-f(s))^2(z-f(t))}\ri)\ud s \ud t.
\eeq
Since, for $M_\si:=\sup_{0\le x\ne y\le 1}\si(x,y)^2$ and for any fixed $z\notin \tS$, the function \[\psi_z:(s,\lam,\lam')\in [0,M_\si+1]\ti \{x\in \R\ste \op{dist}(x,\tS)\le \op{dist}(z,\tS)/2\}^2\longmapsto \f{s}{(z-\lam)^2(z-\lam')}\] is $C\op{dist}(z,\tS)^{-4}$-Lipschitz, for $C$ a universal constant, by Hypothesis (b),
\be\la{estimatedeltan} \delta_n(z)\;=\;O(n^{-1})+ \f{O\lf(\eta_n \ri)}{ \max\{n\eps_n,1\} \op{dist}(z,\tS)^{4} }.\ee

Finally, the expression of $B(z)$ given in \eqref{032016112} implies,
\begin{align}\label{211220161830}
B(z) \leq \f{C}{\op{dist}(z,\tS)^3}
\end{align}

Collecting estimations \eqre{211220161808}, \eqre{estimatedeltan} and \eqre{211220161830}, we conclude.
\epr

\beg{Claim}\la{Claim3}There is  a constant $C$ \st for any   $z\in \C\bck\tS$,
\beq\label{eq:varC(z)}
\E[ |n\eps_n^{-1}C_n(z)|^2]\;\le\;\f{C\eps_n^4}{\op{dist}(z,\tS)^{8}}.
\eeq
\end{Claim}
\bpr
We start by writing for all $z\in\C\bck \tS$\\
$\displaystyle{
\E[ |n\eps_n^{-1}C_n(z)|^2 ] =  \f{\eps_n^{4}}{n^3}  \E\left[ \left| \sum_{i,j,k=1}^n \f{x_{i,j} \ x_{j,k} \ x_{k,i}}{(z-\lambda_n(i))^2 \ (z-\lambda_n(j)) \ (z-\lambda_n(k))} \right|^2 \right]}$
\\
\\
{\color{white}.}\hfill$\displaystyle{ = \f{\eps_n^{4}}{n^3} \sum_{i,j,k,l,m,p=1}^n \f{\mathbb{E}\lf( x_{i,j} \ x_{j,k} \ x_{k,i} \ \ovl{x_{l,m} \ x_{m,p} \ x_{p,l} }\ri)}{(z-\lambda_n(i))^2 \ (z-\lambda_n(j)) \ (z-\lambda_n(k)) \ (\ovl{z}-\lambda_n(l))^2 \ (\ovl{z}-\lambda_n(m)) \ (\ovl{z}-\lambda_n(p))}
}.$

Generically, the set of "edges" $\{(l,m), (m,p), (p,l)\}$ must be equal to the set $\{(i,j), (j,k), (k,i)\}$ in order to get a non zero term. Therefore, the complexity of the previous sum is $O(n^3)$. Note that other non zero terms involving third or fourth moments are much less numerous. Hence,
$$\E[ |n\eps_n^{-1}C_n(z)|^2 ]\leq \f{\eps_n^{4}}{n^3} \times \f{O(n^3)}{\op{dist}(z,\tS)^8}
\; \leq \; \f{C \eps_n^{4}}{\op{dist}(z,\tS)^8}
$$
\epr
\beg{Claim}\la{Claim4}
There is  a constant $C$ \st for any   $z\in \C\bck\R$,
\beq\label{eq:varR(z)}
\E[  
|n\eps_n^{-1}R_n^\eps(z)|^2]\;\le\; \f{ O(n^2 \eps_n^6)}{|\Im(z)|^2 \op{dist}(z,\tS)^8}.
\eeq
\end{Claim}
\bpr
Remind that,
$$R_n^\eps(z):=\f{\eps_n^4}{n}\Tr \ff{z-D} X \ff{z-D} X \ff{z-D}X\ff{z-D}X \ff{z-D^\eps}.$$
Hence,
\begin{align*}
\E[|n\eps_n^{-1}R_n^\eps(z)|^2] &\leq \eps_n^6 \ \E\lf[\lf|\Tr \ff{z-D} X \ff{z-D} X \ff{z-D}X\ff{z-D}X \ff{z-D^\eps}\ri|^2\ri]
\\
&\leq \eps_n^6 \ \E\lf[ \Tr \lf|\lf(\ff{z-D}X\ri)^4\ri|^{2} \times \Tr \lf| \ff{z-D^\eps} \ri|^2 \ri]
\\
&\leq \eps_n^6 \ \E\lf[ \Tr \lf( \lf(\ff{z-D}X\ri)^4 \lf(\ovl{\ff{z-D}X}\ri)^4 \ri) \f{n}{|\text{Im}(z)|^2} \ri]
\\
&\leq \f{n\eps_n^6}{|\text{Im}(z)|^2} \ \E\lf[ \Tr \lf( \lf(\ff{z-D}X\ri)^4 \lf(\ovl{\ff{z-D}X}\ri)^4 \ri) \ri]
\\
&\leq \f{n\eps_n^6}{|\text{Im}(z)|^2} \f{O(n^{5})}{n^4 \ \op{dist}(z,\tS)^8} \; \leq \; \f{O(n^2 \eps_n^6)}{|\text{Im}(z)|^2 \op{dist}(z,\tS)^8}.
\end{align*}
The inequality of the last line takes into account that
\begin{itemize}
\item the $L^8$ norm of the entries of $\sqrt nX$ is uniformly bounded
\item the norm of the entries of $X$ is of order $n^{-1/2}$
\item the norm of the coefficients of $(z-D)^{-1}$ is smaller than $\op{dist}(z,\tS)^{-1}$
\item the complexity of the sum defining the trace is of order $O(n^5)$ since its non-null terms are encoded by four edges trees which have therefore five vertices.
\end{itemize}
\epr

We gather now the results of the previous claims.
\\
\\
For any rate of convergence of $\eps_n$, Claim \ref{Claim1} proves that the process $n\eps_n^{-1}A_n(z)$ converges in distribution to the centered Gaussian variable $Z(z)$. Moreover,
\begin{itemize}
\item if $\eps_n \ll n^{-1}$, then as Claims \ref{Claim2}, \ref{Claim3} and \ref{Claim4} imply that the processes $n\eps_n^{-1}B_n(z)$, $n\eps_n^{-1}C_n(z)$ and $n\eps_n^{-1}R_n^\eps(z)$ converge to $0$ in probability, we can conclude, by Slutsky's theorem, that for any $z\in\C\setminus\R$:
$$ n\eps_n^{-1}\DG_n(z) \; \xrightarrow[n\to\infty]{\text{dist}} \; Z(z) 
$$
\item if $\eps_n \sim \f{c}{n}$, then, as Claims \ref{Claim2}, \ref{Claim3} and \ref{Claim4} imply that the processes $n\eps_n^{-1}B_n(z)$, $n\eps_n^{-1}C_n(z)$ and $n\eps_n^{-1}R_n^\eps(z)$ converge respectively to $cB(z)$, $0$ and $0$ in probability, we can conclude, by Slutsky's theorem, that for any $z\in\C\setminus\R$:
$$n\eps_n^{-1}\DG_n(z) \; \xrightarrow[n\to\infty]{\text{dist}} \; Z(z) + cB(z)$$
\item if $n^{-1}\ll\eps_n\ll n^{-1/3}$, then, as Claims \ref{Claim2}, \ref{Claim3} and \ref{Claim4} imply that the three processes $n\eps_n^{-1}(B_n(z)-\eps_n^2B(z))$, $n\eps_n^{-1}C_n(z)$ and $n\eps_n^{-1}R_n^\eps(z)$ converge to $0$ in probability, we can conclude, by Slutsky's theorem, that for any $z\in\C\setminus\R$:
$$n\eps_n^{-1} \lf(\DG_n(z)-\eps_n^2B(z)\ri) \ninfd Z(z)$$
\end{itemize}

Regarding the convergence in probability \eqref{18mars166}, in the case $n^{-1} \ll \eps_n\ll 1$, Claims \ref{Claim1}, \ref{Claim2}, \ref{Claim3} and \ref{Claim4} imply that the processes $\eps_n^{-2} A_n(z)$, $\eps_n^{-2} B_n(z) - B(z)$, $\eps_n^{-2} C_n(z)$ and $\eps_n^{-2} R_n^\eps(z)$ converge to $0$. 

This finishes the proof of the convergences of Proposition \ref{ConvergenceTheorem}.
\qed

\subsection{Proof of Proposition \ref{StieltjesTheorem}}
Recall that
\[B(z)=\int_{(s,t)\in [0,1]^2}\f{\si^2(s,t)}{(z-f(s))^2(z-f(t))}\ud s\ud t.\]
~\\
Recall that $\rho$  is the density of the push-forward of the uniform measure on $[0,1]$ by the map $f$.
\\
Let $\tau$ be as in Hypotheis (d). We have 
\\
\[B(z) = \ \int_{\mathbb{R}^2} \frac{\tau(s,t) \ \rho(s) \ \rho(t)}{(z-s)^2 \ (z-t)} \ud s \ud t.\]
~\\
By a partial fraction decomposition we have for all $a \neq b$
\[ \frac{1}{(z-a)^2 (z-b)} = \frac{1}{(b-a)^2} \left( \frac{1}{z-b} - \frac{1}{z-a} - \frac{b-a}{(z-a)^2} \right).\]
~\\
Thus, as the Lebesgue measure of the set $\left\{ (y_1,y_2) \in [0,1]^2 \ \ste \ y_1 = y_2 \right\}$ is null, we have
\\
\[ B(z) =  \ \int_{\mathbb{R}^2} \frac{\tau(s,t) \ \rho(s) \ \rho(t)}{(t-s)^2} \left( \frac{1}{z-t} - \frac{1}{z-s} - \frac{t-s}{(z-s)^2} \right) \ud s \ud t.\]
Moreover, for $\varphi_z$ the function $\varphi_z: x \longmapsto \frac{1}{z-x}$, we obtain
\[ B(z) = \ \int_{\mathbb{R}^2} \frac{\tau(s,t) \ \rho(s) \ \rho(t)}{(t-s)^2} \left( \varphi_z(t) - \varphi_z(s) - (t-s) \varphi_z'(s) \right) \ud s \ud t.\]

Now, we want to prove that  $\ds B(z) = -\int_{\mathbb{R}^2} \frac{\tau(s,t) \ \rho(s) \ \rho(t)}{t-s} \ \varphi_z'(s) \ \ud s \ud t$.
 
To do this, we will use a symmetry argument: in fact both terms in $\varphi_z(t)$ and $\varphi_z(s)$ neutralize each other, and it remains only to prove, that we did not remove $\infty$ to $\infty$ and that the remaining term has the desired form.

Let us define 
\[ B^\eta(z):=  \ \int_{\substack{|s-t|>\eta}} \frac{\tau(s,t) \ \rho(s) \ \rho(t)}{(t-s)^2} \left( \varphi_z(t) - \varphi_z(s) - (t-s) \varphi_z'(s) \right) \ud s \ud t.\]
By the Taylor-Lagrange inequality we obtain:
\[
\left| \frac{\tau(s,t) \ \rho(s) \ \rho(t)}{(t-s)^2} \left( \varphi_z(t) - \varphi_z(s) - (t-s) \varphi_z'(s) \right) \right| 
\leq \frac{ \rho(s) \ \rho(t) \ \|\tau(\cdot,\cdot)\|_{L^\infty} \ \|\varphi_z''\|_{L^\infty}}{2}.
\]
So that, since $\rho$ is a density, by dominated convergence, we have
\[ \lim_{\eta \to 0} B^\eta(z) = B(z).\]
Moreover, by symmetry,  for any $\eta$, 
\[ B^\eta(z) =    \int_{\substack{|s-t|>\eta}} \frac{\tau(s,t) \ \rho(s) \ \rho(t)}{t-s} (-\varphi_z'(s) ) \ud s \ud t.\]
So \beqy B(z)&=& \lim_{\eta \to 0}  \int_{|s-t|>\eta} \frac{\tau(s,t) \ \rho(s) \ \rho(t)}{t-s} (-\varphi_z'(s) ) \ud t \ud s\nonumber \\
&=& -\lim_{\eta \to 0}\int_{s\in \R}F_\eta(s) \varphi_z'(s) \ud s\label{1909171}
\eeqy where for $\eta > 0$ and $s \in \R$, we define
\[ F_\eta(s): = \rho(s)\int_{t \in \mathbb{R}\bck[s-\eta,s+\eta]}   \frac{\tau(s,t)  \ \rho(t)}{t-s} \ud t.\]
Note that  that by 
definition of the function $F$ given at \eqre{IntroF}, for any $s$, we have \be\label{1909172}F(s)=\lim_{\eta\to 0}F_\eta(s).\ee
Thus by \eqre{1909171} and \eqre{1909172}, to  conclude the proof of Proposition \ref{StieltjesTheorem}, by dominated convergence, one needs only to state that $F_\eta$ is dominated, uniformly in $\eta$, by an integrable function. 
This follows from the following computation. 

Note first  that by symmetry, we have 
 \be\label{1909173} F_\eta(s) =\rho(s) \int_{t \in \mathbb{R}\bck[s-\eta,s+\eta]}   \frac{\tau(s,t)   \ \rho(t)-\tau(s,s) \  \rho(s)}{t-s} \ud t.\ee 
Let    $M>0$ \st  the  support of the function $\rho$ is contained in $[-M,M]$. Then, for $\eta_0,\al,C$ as in Hypothesis (e),  using the expression of $F_\eta(s)$ given at \eqref{1909173},  we have 
\begin{align*}
\left| F_\eta(s) \right| & \leq \quad 2C\rho(s)\int_{t=s}^{s+\eta_0}|t-s|^{\al-1}\ud t + \int_{t\in [s-2M,s-\eta_0]\cup[s+\eta_0,s+2M]} \left| \frac{  \tau(s,t)\rho(s)\rho(t)}{t-s} \right| \ud t
\\
& \leq \quad \frac{2C\rho(s)}{\al}\eta_0^\al + \frac{1}{\eta_0}\int_{t\in \mathbb{R}} \left| \tau(s,t)\rho(s)\rho(t) \right| \ud t
\\
& \leq \quad \frac{2C\rho(s)}{\al}\eta_0^\al + \frac{\|\tau(\cdot,\cdot)\|_{L^\infty}}{\eta_0}\rho(s).
\end{align*}
\qed

\subsection{A local type convergence result}\label{subsec:local}
One can precise the convergence \eqref{18mars166} by replacing the complex variable $z$ by a complex sequence $(z_n)$ which converges slowly enough to the real axis. This convergence won't be used in the sequel. As it is discussed in \cite{NotesAnttiFlo}, this type of result is a first step towards a local result for the empirical distribution.

\begin{propo}\label{localLemma} Under Hypotheses (a), (b), (f), if $n^{-1}\ll \eps_n\ll 1$, then for any nonreal complex sequence $(z_n)$,   \st
\be\la{71171}\Im(z_n) \gg \max\lf\{ (n\eps_n)^{-1/2} \ , \ \lf(\f{\eta_n}{n\eps_n}\ri)^{1/4} \ , \ \eps_n^{2/5} \ri\} \ee
the following convergence holds \beq\la{18mars165}\eps_n^{-2}\DG_n(z_n) - B(z_n) \ninfP 0 \,.\eeq
\end{propo}

{\bf Remark.} In the classical case where $\displaystyle{\dfrac{\eta_n}{n\eps_n} = \sup_{i\ne j }( |\si_n^2(i,j)-\si^2(i/n,j/n)|+|\lam_n(i)-f(i/n)|)}$ is of order $\dfrac1{n}$, the above assumption boils down to
$\Im(z_n) \gg \max\lf\{ (n\eps_n)^{-1/2} \ , \ \eps_n^{2/5} \ri\}$.

\bpr
Assume $n^{-1} \ll \eps_n\ll 1$. One can directly obtain, for all non-real complex sequences $(z_n)$, that
\begin{itemize}
\item by Claim \ref{Claim1}, if $\op{dist}(z_n,\tS) \gg (n\eps_n)^{-1/2}$, then
\bes\la{Eq2:VarA(z)}\E[ |\eps_n^{-2}A_n(z_n)|^2] \;\leq\; \f{C}{ (n\eps_n)^2 \ \op{dist}(z_n, \tS)^4} \;\ninf \;  0,\ees
\item by Claim \ref{Claim2}, if $\op{dist}(z_n,\tS) \gg \max\lf\{n^{-1/3} \ , \ (\eta_n/(n\eps_n ))^{1/4}\ri\}$, then
$$\mathbb{E}[|\eps_n^{-2}B_n(z_n) - B(z_n)|^2] \;\le \;\f{C}{n^{2}\op{dist}(z_n, \tS)^6}+\f{C \eta_n^2 }{(n\eps_n)^{2} \ \op{dist}(z_n, \tS)^8} \ninfe 0,$$
\item by Claim \ref{Claim3}, if $\op{dist}(z_n,\tS) \gg \lf(\eps_n/n\ri)^{1/4}$, then
$$ 
\E[|\eps_n^{-2}C_n(z_n)|^2] \;\le\;\f{C\eps_n^2}{n^{2} \ \op{dist}(z_n,\tS)^{8}} \;  \ninfe\;  0,
$$
\item by Claim \ref{Claim4}, if $ |\Im(z_n)| \op{dist}(z_n,\tS)^4 \gg \eps_n^{2}$, then
$$
\E[|\eps_n^{-2}R_n^\eps(z_n)|^2] \; \leq \f{O(\eps_n^4)}{|\Im(z_n)|^2 \op{dist}(z_n,\tS)^8} \; \ninfe\;  0. $$
\end{itemize}
Therefore, when
$$\op{dist}(z_n,\tS) \gg \max\lf\{ (n\eps_n)^{-1/2} \ , \ n^{-1/3} \ , \ \lf(\f{\eta_n}{n\eps_n}\ri)^{1/4} \ , \ \lf(\f{\eps_n}{n}\ri)^{1/4}  \ri\} \ \text{ and } \ |\Im(z_n)| \op{dist}(z_n,\tS)^4 \gg \eps_n^{2},$$
the four processes, $\eps_n^{-2}A_n(z_n)$, $\eps_n^{-2}B_n(z_n) - B(z_n)$, $\eps_n^{-2}C_n(z_n)$ and $\eps_n^{-2}R_n^\eps(z_n)$ converge to $0$ in probability.
Since $\op{dist}(z_n,\tS) \geq\Im(z_n)$, the above condition is implied by
\[\Im(z_n) \gg \max\lf\{ (n\eps_n)^{-1/2} \ , \ n^{-1/3} \ , \ \lf(\f{\eta_n}{n\eps_n}\ri)^{1/4} \ , \ \lf(\f{\eps_n}{n}\ri)^{1/4}, \ \eps_n^{2/5} \ri\}.\]
Observing finally that the two terms $n^{-1/3}$ and $\lf(\f{\eps_n}{n}\ri)^{1/4}$ are dominated by the maximum of the three other ones, we conclude the proof.
\epr

\subsection{Possible extensions to larger  $\eps_n$}\label{subsec:extensions}

The convergence in distribution result of Theorem \ref{ConvergenceTheorem00} is valid for $\eps_n\ll n^{-1/3}$ but fails above $n^{-1/3}$. Let us consider, for example, the case where $n^{-1/3}\ll\eps_n\ll n^{-1/5}$. In this case, the contribution of the first term $A_n(z)$ in the expansion of  $\DG_n(z)$ which yields the random limiting quantity, is dominated {\it not only} by the term $B_n(z)$ as it used to be previously. It is also dominated by a further and smaller term $D_n(z)$ of the expansion 
 \[\DG_n(z) = A_n(z)+B_n(z)+C_n(z)+D_n(z)+E_n(z)+R_n^\eps,\]
  with: 
  \beq A_n(z)&:=&  \f{\eps_n}{n} \Tr\ff{z-D}X\ff{z-D}\\
&\vdots&
  \\
  E_n(z)&:=&\f{\eps_n^5}{n}\Tr \ff{z-D}X\ff{z-D}X\ff{z-D}X\ff{z-D}X\ff{z-D}X\ff{z-D}\\
  R_n^\eps(z)&:=&\f{\eps_n^6}{n}\Tr \ff{z-D}X \ff{z-D}X \ff{z-D}X \ff{z-D}X\ff{z-D}X\ff{z-D}X\ff{z-D^\eps}.\eeq
  
In this case, the random term $Z(z)$ is still produced by $A_n(z)$ and has an order of magnitude of $\eps_n/n$. Meanwhile, the term $D_n(z)$ writes
$$ D_n(z) := \f{\eps_n^4}{n^3} \sum_{i,j,k,l=1}^n \f{x_{i,j} \ x_{j,k} \ x_{k,l} \ x_{l,i}}{ (z-\lambda_n(i))^2 \ (z-\lambda_n(j)) \ (z-\lambda_n(k)) \ (z-\lambda_n(l)) }. $$ 
All the indices satisfying $j=l$ contribute to the previous sum, since they produce a term in $|x_{i,l}|^2 |x_{k,l}|^2$. Their cardinality is of order $n^3$. Therefore, the  term $D_n(z)$ is of order $\eps_n^4$ which prevails on the order $\eps_n/n$ of $A_n(z)$, as soon as $\eps_n\gg n^{-1/3}$.
One can also observe that the odd terms $C_n(z)$ and $E_n(z)$ in the expansion are negligible with respect to $A_n(z)$ due to the fact that  the entries $x_{i,j}$ are centered. 
One can then state an analogous result to Proposition \ref{ConvergenceTheorem}, but the deterministic limiting term $D(z)$ arising from $D_n(z)$ does not find a nice expression as the image of $\varphi_z$ by a linear form as it was the case  for $B(z)$ in Proposition  \ref{StieltjesTheorem}. Therefore we did not state an extension of Theorem \ref{ConvergenceTheorem00}. 

More generally, for all positive integer $p$, when $n^{-1/(2p-1)}\ll \eps_n\ll n^{-1/(2p+1)}$, the expansion will contain $p$ deterministic terms, produced by the even variables, $B_n(z)$, $D_n(z)$, $F_n(z)$, $H_n(z)$ $\dots$  All the other odd terms, $C_n(z)$, $E_n(z)$, $G_n(z)$ $\dots$ being negligible due to the centering of the entries. The limits of the even terms $B_n(z)$, $D_n(z)$, $F_n(z)$, $H_n(z)$ $\dots$ can be expressed thanks to operator-valued free probability theory, using the results of \cite{Sh96} (namely, Th. 4.1), but expressing these limits as the images of $\varphi_z$ by  linear forms is a quite involved combinatorial problem that we did not solve yet.

\section{Convergence in probability in the semi-perturbative regime}\la{sec:probability}

Our goal now is to extend the convergence in probability result \eqref{18mars166} of Proposition \ref{ConvergenceTheorem}, proved for test functions $\varphi_z(x) := \ff{z-x}$, to any $\mathcal{C}^6$ and compactly supported function on $\mathbb{R}$. We do it in the following lemma by using the Helffer-Sj\"ostrand formula which is stated in Proposition \ref{prop:HS} of the Appendix.

\begin{lem}\label{lem:proba}
If $n^{-1} \ll \eps_n \ll 1$, then, for any compactly supported $\mathcal{C}^6$ function $\phi$ on $\R$,
\[ \eps_n^{-2}(\mu_n^\eps - \mu_n)(\phi) \xrightarrow[n\to\infty]{P} -\int\phi'(s) F(s) \ \ud s\, . \]
\end{lem}

\bpr Let us introduce the Banach space $\mc{C}^1_{\op{b},\op{b}}$ of bounded  $\mc{C}^1$ functions on $\R$ with  bounded derivative, endowed with the  norm $\|\phi\|_{\mc{C}^1_{\op{b},\op{b}}}:=\|\phi\|_\infty+\|\phi'\|_\infty$.

On this space, let us define the random continuous   linear form
\bes\label{30nov1554} \Pi_n(\phi):= \eps_n^{-2}(\mu_n^\eps - \mu_n)(\phi) + \int\phi'(s) F(s) \ \ud s. \ees 

Convergence \eqref{18mars166} of Proposition \ref{ConvergenceTheorem} can now be formulated as 
$$\forall z\in\C\setminus\R, \qquad \Pi_n(\varphi_z) \xrightarrow[n\to\infty]{P} 0.$$

Actually, we can be more precise by adding the upper bounds of Claims \ref{Claim1}, \ref{Claim2}, \ref{Claim3} and \ref{Claim4}, and obtain, uniformly in $z$,
\begin{align}
\E[|\Pi_n(\varphi_z)|^2] \; &= \; \E[| \eps_n^{-2}\DG_n(z)-B(z)|^2] \nonumber
\\
\; &\leq \; \f{(n\eps_n)^{-2}}{\min\lf( \op{dist}(z,\tS)^4 \ , \ \op{dist}(z,\tS)^8 \ , \ |\Im(z)|^2 \op{dist}(z,\tS)^8 \ri)}. \label{1504091216}
\end{align}
\\
Now, let  $\phi$ be  a compactly supported $\mc{C}^{6}$ function on $\R$ and let us introduce the almost analytic extension of degree $5$ of $\phi$ defined by \begin{equation*}
\forall z=x+\ii y\in\C, \qquad \widetilde \phi_5(z) \ \deq \ \sum_{k = 0}^5 \frac{1}{k !} (\ii y)^k \phi^{(k)}(x)\,.
\end{equation*}
An elementary computation gives, by successive cancellations, that
\be\label{eq23111614112} \bar \partial \widetilde \phi_5(z) = \ff{2}\lf( \partial_x  + \ii \partial_y \ri) \widetilde \phi_5(x+\ii y) = \ff{2\times 5!} (\ii y)^5 \phi^{(6)}(x).\ee

Furthermore, by Helffer-Sj\"ostrand formula (Proposition \re{prop:HS}), for $\chi \in \mathcal{C}^\infty_c(\C;[0,1])$   a smooth cutoff function with value one on the support of $\phi$,
\begin{equation*}
\phi(\cdot) \;=\; -\frac{1}{\pi} \int_{\C} \frac{\bar \partial (\widetilde \phi_5(z) \chi(z))}{y^5}y^5\varphi_z(\cdot) \, \dd^2 z\,
\end{equation*}
where $\dd^2 z$ denotes the Lebesgue measure on $\C$.

Note that by \eqref{eq23111614112}, $z \mapsto \one_{y\ne 0} \frac{\bar \partial (\widetilde \phi_5(z) \chi(z))}{y^5}$ is a continuous compactly supported function and that $z\in \C\mapsto  \one_{y\ne 0} y^5\vfi_z\in \mc{C}^1_{\op{b},\op{b}}$ is continuous, hence,
\[ \Pi_n(\phi) = \frac{1}{\pi} \int_{\C}  \f{\bar \partial (\widetilde \phi_5(z) \chi(z))}{y^5} \ y^5\Pi_n(\varphi_z) \, \dd^2 z. \]
Therefore,
\begin{align*}
\E\lf( \lf| \Pi_n(\phi) \ri|^2 \ri) &= \E\lf( \lf| \frac{1}{\pi} \int_{\C}  \f{\bar \partial (\widetilde \phi_5(z) \chi(z))}{y^5} \ y^5\Pi_n(\varphi_z) \, \dd^2 z \ri|^2  \ri)
\\
& \leq \E\lf( \frac{1}{\pi^2} \int_{\C} \lf| \f{\bar \partial (\widetilde \phi_5(z) \chi(z))}{y^5} \ y^5\Pi_n(\varphi_z) \ri|^2 \, \dd^2 z  \ri)
\\
&= \frac{1}{\pi^2} \int_{\C} \lf| \f{\bar \partial (\widetilde \phi_5(z) \chi(z))}{y^5}\ri|^2 \ y^{10} \ \E\lf(\lf| \Pi_n(\varphi_z)\ri|^2\ri)  \, \dd^2 z\, .
\end{align*}

Since the function $\lf| \f{\bar \partial (\widetilde \phi_5(z) \chi(z))}{y^5}\ri|^2$ is continuous and compactly supported and that, by \eqre{1504091216}, for $n^{-1} \ll \eps_n \ll 1$, uniformly in $z$,
\[ y^{10} \ \E\lf(\lf| \Pi_n(\varphi_z)\ri|^2\ri) \leq y^{10} \ \f{o(1)}{\min(y^4, \ y^{10})} \ninfe 0. \]

Thus, for any compactly supported $\mc{C}^6$ function on $\R$,
\[ \E\lf( \lf| \Pi_n(\phi) \ri|^2 \ri) \leq \frac{1}{\pi^2} \int_{\C} \lf| \f{\bar \partial (\widetilde \phi_5(z) \chi(z))}{y^5}\ri|^2 \ y^{10} \ \E\lf(\lf| \Pi_n(\varphi_z)\ri|^2\ri)  \, \dd^2 z \ninfe 0 \]
which implies that $\Pi_n(\phi)$ converges to $0$ in probability.
\epr

\section{Convergence in distribution towards the Gaussian variable $Z_\phi$}\la{sec:clt} 

The purpose of this section is to extend the convergences in distribution of Proposition \ref{ConvergenceTheorem}, from test functions of the type $\varphi_z:= \ff{z-x}$, to compactly supported  $\mathcal{C}^6$ functions on $\mathbb{R}$. To do so,  we will use an extension lemma of Shcherbina and Tirozzi, stated in Lemma \ref{extensionLemma} of the Appendix, which concerns the convergence of a sequence of {\it centered} random fields with uniformly bounded variance. Hence, we need to show first that our non centered random sequence is not far from being centered, which is done in subsection \ref{coincidence} by using again the Helffer-Sj\"ostrand formula (\ref{prop:HS}). In subsection \ref{subsec:shcherbina}, we dominate the variance of this centered random field thanks to another result of Shcherbina and Tirozzi stated in Proposition \ref{extensionLemma2} of the Appendix. Subsection \ref{subsec:collect} collects the preceding results to conclude the proof.

\subsection{Coincidence of the expectation of $\mu_n^\eps$ with its deterministic approximation}\label{coincidence}

The asymptotic coincidence of the expectation of $\mu_n^\eps$ with its deterministic approximation is the content of next lemma:

 \beg{lem}\la{lem1731616h}Let us define, for   $\phi$   a $\mc{C}^1$ function on $\R$,
$$\La_n(\phi):=\beg{cases}n\eps_n^{-1}\lf(\E[\mu_n^\eps(\phi)]-\mu_n(\phi)\ri)&\trm{ if $\eps_n\ll n^{-1}$,}\\ \\ 
n\eps_n^{-1}\lf(\E[\mu_n^\eps(\phi)]-\mu_n(\phi)+\eps_n^2\int \phi'(s)F(s)\ud s\ri)&\trm{if } \eps_n\sim c/n \trm{ or } n^{-1} \ll \eps_n \ll n^{-1/3}\, .
\end{cases}$$
  Then, as $n\to\infty$, for any compactly supported $\mc{C}^6$ function $\phi$ or any $\phi$ of the type $\vfi_z(x)=\ff{z-x}$, $z\in \C\bck \R$, we have   
$$ \La_n(\phi)\ninf 0.$$
\en{lem}

\bpr
First note that, as the variables $x_{i,j}$ are centered, $\E[A_n(z)] = 0$. Moreover, by adding the renormalized upper bounds of Claims \ref{Claim2}, \ref{Claim3} and \ref{Claim4} one can directly obtain the two following inequalities for any $z\in\C\setminus\R$:
\bgt
\ite If $ \eps_n\ll n^{-1}$, then 

	\beq \label{ineq:esp} |\Lambda_n(\varphi_z)| &=& n\eps_n^{-1}|
	\E[\DG_n(z)]| 
	\\ \\ &\leq & n\eps_n^{-1}\lf( |\E[A_n(z)]| + \E[|B_n(z)|] + \E[|C_n(z)|] + \E[|R_n^\eps(z)|] \ri)
	\\ \\ &\leq & \f{ C(n\eps_n + \eta_n) }{ \min\left\{ \op{dist}(z, \tS)^3, \op{dist}(z, \tS)^4,|\Im(z)| \op{dist}(z, \tS)^4\right\} } \ninfe 0\,.
	 \eeq
	 
\ite If $\eps_n\sim c/n$ or $n^{-1}\ll \eps_n\ll n^{-1/3}$, then

	\beq \label{ineq:esp2} |\Lambda_n(\varphi_z)| &=& n\eps_n^{-1}|\E[\DG_n(z)-\eps_n^{2}B(z)]|
	\\ \\ &\leq & n\eps_n^{-1}\lf( |\E[A_n(z)]| + \E[|B_n(z)-\eps_n^{2}B(z)|] + \E[|C_n(z)|] + \E[|R_n^\eps(z)|] \ri)
	\\ \\ &\leq & \f{ C(\eps_n + \eta_n + n\eps_n^3) }{ \min\left\{ \op{dist}(z, \tS)^3, \op{dist}(z, \tS)^4,|\Im(z)| \op{dist}(z, \tS)^4\right\} } \ninfe 0 \, .
	 \eeq
\ent

Hence, in all cases, $ \La_n(\varphi_z)\ninf 0 $.

The extension of this result to compactly supported $\Cc^6$ test functions on $\R$ goes the same way as for $\Pi_n$ in the proof of Lemma \ref{lem:proba}. 
\epr

\subsection{Domination of the variance of $\mu_n^\eps$}\label{subsec:shcherbina}
The second ingredient goes through a domination of the variance of $\mu_n^\eps(\phi)$:
\beg{lem}\la{lemmamajorationvariance}Let $s>5$. There is a constant $C$ \st  for each $n$ and each $\phi\in\Hc_s$, 
$$\Var(n\eps_n^{-1}\mu_n^\eps( \phi))\;\le\; C\|\phi\|_{\Hc_s}^2\,.$$
\en{lem}

\bpr By Proposition \re{extensionLemma2}, it suffices to prove that
$$\int_{y=0}^{\infty}y^{2s-1}\me^{-y}\int_{x\in \R}\Var(\eps_n^{-1} \Tr ((x+\ii  y-D_n^\eps)^{-1}))\ud x\ud y$$
are bounded independently of $n$. 

Note that for $\DG_n(z)$ defined in \eqre{032016113},
$$\Var(\eps_n^{-1} \Tr ((z-D_n^\eps)^{-1}))=n^2\eps_n^{-2} \Var(\DG_n(z)).$$
Moreover, the sum of the inequalities of Claims \ref{Claim1}, \ref{Claim2}, \ref{Claim3} and \ref{Claim4} yields
\beq
\Var( n\eps_n^{-1}\DG_n(z)) \; &\le& \; \f{C}{\op{dist}(z, \tS)^4} + \f{C}{|\Im(z)|^2 \op{dist}(z, \tS)^8} .
\eeq

Let $M>0$ \st $\tS\subset [-M,M]$. Then $$
 \op{dist}(z,\tS)\ge \beg{cases} y&\trm{ if $|x|\le M$}, \\
 \sqrt{y^2+(|x|-M)^2}&\trm{ if $|x|> M$}.
 \end{cases}$$
Thus $ {\op{dist}(z,\tS)}\ge y$ if $|x|\le M$ and, for $|x|> M$, 
$$ 
\ff{\op{dist}(z,\tS)} \le \f{y^{-1}}{\sqrt{ 1+((|x|-M)/y)^2 }}
$$
and for any $y>0$,
\begin{align*}
\int_{x\in \R}\Var(n\eps_n^{-1}\DG_n(x+\ii y))\ud x \; &\le \; 2CM (y^{-10} + y^{-4}) + 2C \int_0^{+\infty}\f{y^{-4}}{(1+(\f{x}{y})^2)^{2}} + \f{y^{-10}}{(1+(\f{x}{y})^2)^{4}} \ud x
\\
&\le \; 2CM (y^{-10} + y^{-4}) + C \lf( \f{\pi}{2} y^{-3} + \f{5\pi}{16} y^{-9}\ri)
\\
&\le \; k \lf(y^{-10}+  y^{-3}\ri),
\end{align*}
for a suitable constant $k$.

We deduce that, as soon as $2s-10>0$, \ie $s>5$, 
\[ \int_{y=0}^{\infty}y^{2s-1}\me^{-y}\int_{x\in \R}\Var(\eps_n^{-1} \Tr ((x+\ii  y-D_n^\eps)^{-1})) \ud x\ud y \ \le \  k \int_{0}^{\infty}y^{2s-1}e^{-y} (y^{-10}+y^{-3}) \ud y \ < \ \infty. \]

\epr

\subsection{Proof of the convergences in distribution of Theorem \ref{ConvergenceTheorem00}}\label{subsec:collect}

Since we have proved in Lemma \ref{lem1731616h} that for all compactly supported $\Cc^6$ function  $\phi$, the deterministic term $\mu_n(\phi)$ could be replaced by $\E[\mu_n^\eps(\phi)]$, we only have to prove, that for all $\phi\in \Cc^6$,
$$n\eps_n^{-1}(\mu_n^\eps(\phi)-\E[\mu_n^\eps(\phi)])  \ninfd  Z_\phi.$$  
 For the time being, we know this result to be valid for functions $\phi$  belonging to the space $\Lc_1$, defined as the linear span of the family of functions $\varphi_z(x):=\ff{z-x}$, $z\in \C\bck\R$.

By applying Lemma \ref{extensionLemma} to the centered field $\mu_n^\eps-\E[\mu_n^\eps]$, we are going to extend the result from the space $\Lc_1$ to the Sobolev space $(\Hc_s, \|\cdot\|_{\Hc_s})$ with $s\in (5,6)$. Note that, since $s<6$, this latter space contains the space of $\Cc^6$ compactly supported functions (see \cite[Sec. 7.9]{Hormander}).

It remains to check the two hypotheses of Lemma \ref{extensionLemma}. First, the subspace $\Lc_1$ is dense in every space  $(\Hc_s, \|\cdot\|_{\Hc_s})$.  This is the content of Lemma \re{113164} of the Appendix.  Second, by Lemma \re{lemmamajorationvariance}, since $s>5$, $\Var(n\eps_n^{-1}\mu_n^\eps( \phi))\le C\|\phi\|_{\Hc_s}^2$ for a certain constant $C$.

This concludes the proof.

\section{Appendix}

The reader can find here the results we use along the paper, namely the Helffer-Sj\"ostrand formula, the CLT extension lemma of Shcherbina and Tirozzi and a functional density lemma with its proof.

\subsection{Helffer-Sj\"ostrand formula}
The proof of the following formula can be found, e.g. in \cite{NotesAnttiFlo}. 
\begin{propo}[Helffer-Sj\"ostrand formula] \label{prop:HS}
Let $n \in \mathbb{N}$ and $\phi \in  \mathcal{C}^{p+1}(\mathbb{R})$. We define the \emph{almost analytic extension of $\phi$ of degree $p$} through
\begin{equation*}
\widetilde \phi_p(x + \ii y) \ \deq \ \sum_{k = 0}^p \frac{1}{k !} (\ii y)^k \phi^{(k)}(x)\,.
\end{equation*}
Let $\chi \in \mathcal{C}^\infty_c(\C;[0,1])$ be a smooth cutoff function. Then for any $\lambda \in \R$ satisfying $\chi(\lambda) = 1$ we have
\begin{equation*}
\phi(\lambda) \;=\; \frac{1}{\pi} \int_\C \frac{\bar \partial (\widetilde \phi_p(z) \chi(z))}{\lambda - z} \, \dd^2 z\,,
\end{equation*}
\\
where $\dd^2 z$ denotes the Lebesgue measure on $\mathbb{C}$ and $\bar \partial \deq \frac{1}{2} (\partial_x + \ii \partial_y)$ is the antiholomorphic derivative.
\end{propo}

\subsection{CLT extension lemma}\la{sec:extensionlemmas}

The following CLT extension lemma is borrowed from the paper of Shcherbina and Tirozzi \cite{tirozzi}. We state here the version that can be found in the Appendix of \cite{ACFTCL}.

\beg{lem}\label{extensionLemma}
Let $(\Lc, \|\,\cdot\,\|)$ be a normed space with a dense subspace $\Lc_1$ and, for each $n\ge 1$, $(N_n(\phi))_{\phi\in \Lc}$ a collection of real   random variables such that: \bgt\ite for each $n$ , $\phi\longmapsto N_n(\phi)$ is linear, 
 \ite  for each $n$ and each $\phi\in \Lc$, $ \E[N_n(\phi)]=0 $,
 \ite there is a constant $C$ \st  for each $n$ and each $\phi\in\Lc$, $\Var(N_n(\phi))\le C\|\phi\|^2$, 
 \ite there is a quadratic form $V:\Lc_1\to\R_+$ \st for any 
 $\phi\in \Lc_1$, we have the convergence in distribution $N_n(\phi)\ninf \mc{N}(0, V(\phi))$.
 \ent
 Then $V$ is continuous on $\Lc_1$, can  (uniquely) be continuously extended to $\Lc$ and   for any $\phi\in \Lc$, we have the convergence in distribution $N_n(\phi)\ninf \mc{N}(0, V(\phi))$.
\en{lem}

One of the assumptions of previous lemma concerns a variance domination.  The next proposition provides a tool in order to check it. 
Let us first remind  the definition of the Sobolev space $\Hc_s$. 
For $\phi\in L^1(\R,\ud x)$, we define $$\hfi(k):= \int \me^{ \ii  kx}\phi(x)\ud x\qquad\qquad\trm{ ($k\in \R$)}$$ and, for $s>0$, $$\|\phi\|_{\Hc_s}:=\|k\longmapsto (1+2|k|)^s\,\hfi(k)\|_{L^2}.$$
We define the Sobolev space  $\Hc_s$ as  the  set of functions  with finite $\|\cdot\|_{\Hc_s}$ norm. Let us now state Proposition 2 of the paper \cite{tirozzi2} of Shcherbina and Tirozzi.

\beg{propo}\label{extensionLemma2} For any $s>0$, there is a constant $C=C(s)$ \st for any $n$, any $n\ti n$ Hermitian random matrix  $M$, 
and any $\phi\in \Hc_s$, 
we have \bes\la{9121311}\Var(\Tr \phi(M))\;\le \;C \|\phi\|_{\Hc_s}^2\int_{y=0}^{\infty}y^{2s-1}\me^{-y}\int_{x\in \R}\Var(\Tr ((x+\ii  y-M)^{-1}))\ud x\ud y.
\ees
\en{propo}

\subsection{A density lemma}
We did not find Lemma \re{113164} in the literature, so we provide its proof.
Recall that for any $z\in \C\bck\R$, $$\varphi_z(x)=\ff{z-x}.$$
\beg{lem}\la{113165}For any $z\in \C\bck\R$, we have, in the $L^2$ sense, \be\la{eq:113165}\wt{\varphi_z}\;=\;(t\longmapsto -\op{sgn}(\Im z)2\pi\ii\one_{\Im(z)t>0}\me^{\ii tz})\ee and $\varphi_z$ belongs to each $\Hc_s$ for any  $s\in \R$.
\en{lem}

\bpr 
  It is well known that if $\Re z>0$, then $\ds \f{1}{z}= \int_{t=0}^{+\infty} \me^{-tz}\ud t.$  
  
  Let   $z=E+\ii\eta$, $E\in \R, \eta>0$. For any $\xi\in \R$, we have 
  \bes\la{68123bis} \varphi_z(\xi)=\f{-\ii}{\ii(\xi-z)}= -\ii\int_{t=0}^{+\infty} \me^{-\ii t(\xi-z)}\ud t= -\ii\int_{t=0}^{+\infty} \me^{-\ii t\xi} \me^{\ii tz}\ud t   .\ees We deduce \eqre{eq:113165} for $\Im z>0$. The general result can be deduced by complex conjugation.  
\epr

  \beg{lem}\la{113164}Let $\Lc_1$ denote the linear span of the functions $\varphi_z(x):=\ff{z-x}$, for $z\in \C\bck\R$.  Then the space $\Lc_1$ is dense in $\Hc_s$ for any $s\in \R$.
  \end{lem}

  \bpr   We know, by Lemma \re{113165}, that $\Lc_1\subset \Hc_s$. 
  Recall first the definition of the Poisson kernel, for $E\in \R$ and $\eta>0$, $$P_\eta(E)=\ff{\pi}\f{\eta}{E^2+\eta^2}=\ff{2\ii\pi}\lf(\vfi_{\ii\eta}(E)-\vfi_{-\ii\eta}(E)\ri)$$ and that,   by Lemma \re{113165}, $$\wt{P_\eta}(t)=\me^{-\eta|t|}.$$ Hence for any $f\in \Hc_s$, we have $$\|f-P_\eta*f\|_{\Hc_s}^2=\int (1+2|x|)^{2s}|\wt{f}(x)|^2(1-\me^{-\eta|x|})^2\ud x,$$ so that, by dominated convergence, $P_\eta*f\lto f$ in $\Hc_s$ as $\eta\to 0$. 
  
To prove Lemma \re{113164}, it suffices to prove that any smooth compactly supported function can be approximated, in $\Hc_s$, by   functions of $\Lc_1$. So let $f$ be a smooth compactly supported function. By what precedes, it suffices to prove that for any fixed $\eta>0$, $P_\eta*f$ can be approximated, in $\Hc_s$, by   functions of $\Lc_1$.
For $x\in \R$, \beq P_\eta*f(x)&=&\ff{\pi}\int f(t)\f{\eta}{\eta^2+(x-t)^2}\ud t
  \\
 &=&-\ff{\pi}\int f(t)\Im(\vfi_{t+\ii\eta}(x))\ud t  \\
 &=&\ff{2\pi\ii}\int f(t)(\vfi_{t-\ii\eta}(x)-\vfi_{t+\ii\eta}(x))\ud t.
  \eeq
Without loss of generality, one can suppose that the support of $f$ is contained in $[0,1]$.
 Then, for any $n\ge 1$,
 \beqy\la{Mogwai1} P_\eta*f(x) &=&\ff{2n\pi\ii}\sum_{k=1}^n f(\f{k}{n})\lf(\vfi_{\f{k}{n}-\ii\eta}(x)-\vfi_{\f{k}{n}+\ii\eta}(x)\ri)+R_n(x) 
  \eeqy
  where for $[t]_n:=\lceil nt\rceil/n$,  $$R_n(x)= \ff{2\pi\ii}\int f(t)\lf(\vfi_{t-\ii\eta}(x)-\vfi_{t+\ii\eta}(x)\ri)-f([t]_n)\lf(\vfi_{[t]_n-\ii\eta}(x)-\vfi_{[t]_n+\ii\eta}(x)\ri)\ud t.$$
  The error term $R_n(x)$ rewrites \beq R_n(x)&=& \ff{2\pi\ii}\int (f(t)-f([t]_n))(\vfi_{t-\ii\eta}-\vfi_{t+\ii\eta})(x)\ud t\\
  &&+\ff{2\pi\ii}\int f([t]_n)(\vfi_{t-\ii\eta}-\vfi_{[t]_n-\ii\eta}+\vfi_{t+\ii\eta}-\vfi_{[t]_n+\ii\eta})(x)\ud t.\eeq
  
  Now, note that for any $t\in \R$ and $\eta\in \R\bck\{0\}$, we have 
  by Lemma \re{113165},  $$\wt{\vfi_{t+i\eta}}=(x\mapsto - \op{sgn}(\eta)2\pi\ii \, \one_{\eta x>0} \, \me^{\ii xz}),$$ so that when, for example, $\eta>0$, for any $t\in \R$,  $$\|\vfi_{t+\ii\eta}\|_{\Hc_s}^2=4\pi^2\int_0^\infty (1+2|x|)^{2s}\me^{-2\eta x}\ud x$$ does not depend on $t$ and for any $t,t'\in \R$,  \beq \|\vfi_{t+\ii\eta}-\vfi_{t'+\ii\eta}\|_{\Hc_s}^2&=&4\pi^2\int_0^\infty (1+2|x|)^{2s}|\me^{\ii tx}-\me^{\ii t' x}|^2\me^{-2\eta x}\ud x\\ &=&4\pi^2\int_0^\infty (1+2|x|)^{2s}|\me^{\ii (t-t')x}-1|^2\me^{-2\eta x}\ud x\eeq depends only on $t'-t$ end tends to zero (by dominated convergence) when $t'-t\to 0$.
  
  We deduce that $\|R_n\|_{\Hc_s}\lto 0$ as $n\to\infty$, which closes the proof, by \eqre{Mogwai1}. 
  \epr

\vspace{.8cm}
\noindent{\bf Acknowledgements.} We  thank  Jean-Philippe Bouchaud,  Guy David  and Vincent Vargas  for some fruitful discussions. We are also glad to thank the GDR MEGA for partial support.

\end{document}